\documentclass[a4paper,twoside,12pt]{amsart}

\usepackage{amsmath, amsfonts, amssymb}
\usepackage{graphicx}
\usepackage{amsthm}
\usepackage{amsmath}
\usepackage[ansinew]{inputenc}
\usepackage{xy}

\usepackage{amssymb}

\usepackage{epsfig}

\newtheorem{The}{Theorem}[section]
\newtheorem{Lem}{Lemma}[section]

\newtheorem{Cor}[The]{Corollary}
\newtheorem{Pro}[Lem]{Proposition}

\newtheorem{Rem}[Lem]{Remark}
\newtheorem{Defi}[Lem]{Definition}
\newtheorem{Exam}[Lem]{Example}
\newtheorem{Not}[Lem]{Notation}

\def\a{\alpha}

\def\b{\beta}
\def\C{\mathbf C}
\def\D{\Delta}

\def\f{\phi}
\def\g{\gamma}

\def\l{\lambda}

\def\O{\Omega}
\def\p{\pi}

\def\R{\mathbf R}
\def\r{\rho}

\def\s{\sigma}
\def\t{\tau}
\def\tr{{\underline{t}}}

\def\w{\omega}
\def\Z{\mathbf Z}

\title{ Multi-Harnack smoothings of real plane branches}
\author{P.D. Gonz\'alez P\'erez}

\address{Departamento de Algebra. Facultad de Ciencias Matem\'aticas. Universidad Complutense de Madrid.
Plaza de las Ciencias 3. 28040. Madrid. Spain.}

\email{pgonzalez@mat.ucm.es}
\thanks{Gonz\'alez P\'erez is supported by {\em Programa Ram\'on y Cajal} and
MTM2007-6798-C02-02 grants of {\em Ministerio de Educaci\'on y
Ciencia}, Spain.}

\author{J.-J. Risler}
\address{Institut de Mathématiques de Jussieu. \'Equipe Analyse Alg\'ebrique. Case 247, 4 Place Jussieu, 75252 Paris Cedex}

\email{risler@math.jussieu.fr}



\begin{document}

\maketitle

\section*{Introduction}

The 16th problem of Hilbert addresses the determination and the
understanding of the  possible topological types of smooth real
algebraic curves of a given degree in the projective plane $\R
P^2$. This paper is concerned with a local version of this
problem: given a germ $(C, 0)$ of real algebraic plane curve
singularity, determine the possible topological types of the
smoothings of $C$. A {\em smoothing} of $C$ is a real analytic
family $C_t \subset B$, for $t \in [0,1]$, such that $C_0 = C$ and
$C_t$ is non singular and transversal to the boundary of a Milnor
ball $B$ of the singularity $(C,0)$  for $0 < t \ll 1 $. In this
case the real part $\R C_t$ of $C_t$ consists of finitely many
ovals and non closed components in the Milnor ball.

In the algebraic case it was shown by Harnack that a real
projective curve of degree $d$ has at most $\frac{1}{2}(d-1)(d-2)
+1$ connected components. A curve with this number of components
is called a $M$-curve. In the local case there is a similar bound,
depending on the number of real branches of the singularity (see
Section \ref{smoothing}), which arises from the application of the
classical topological theory of Smith. A smoothing which reaches
this bound on the number of connected components is called a
$M$-smoothing. It should be noticed that in the local case
$M$-smoothings do not always exists (see \cite{KOS}). One relevant
open problem in the theory is to determine the actual maximal
number of components of a smoothing of $(C, 0)$, for $C$ running
in a suitable form of equisingularity class refining the classical
notion of Zariski of equisingularity class in the complex world
(see \cite{KRS}).

Quite recently Mikhalkin has proved a beautiful topological
rigidity property of those  $M$-curves in $\R P^2$ which are
embedded in {\em maximal position} with respect to the coordinate
lines (see \cite{Mikhalkin}). His result, which holds more
generally, for those $M$-curves in projective toric surfaces which
are cyclically in maximal position with respect to the toric
coordinate lines, is proved by analyzing the topological
properties of the associated amoebas. The {\em amoeba} of a curve
$C$ is the image of the points $(x, y) \in (\C^*)^2$ in the curve
by the map $Log: (\C^*)^2 \rightarrow \R^2$, given by $(x, y)
\mapsto (\log |x|, \log |y|)$. Conceptually, the amoebas are
intermediate objects which lay in between classical algebraic
curves and tropical curves. See \cite{FPT}, \cite{GKZ},
\cite{Mikhalkin}, \cite{MR}, \cite{PR} and \cite{I} for more on
this notion and its applications.

In this paper we study smoothings of a {\em real plane branch}
singularity $(C,0)$, i.e., the germ $(C, 0)$ is analytically
irreducible in $(\C^2,0)$ and admits a real Newton-Puiseux
parametrization. Risler proved that any such germ $(C,0)$ admits a
$M$-smoothing with the maximal number   ovals, namely $
\frac{1}{2}\mu (C)_0$, where $\mu$ denotes the Milnor number. The
technique used, called nowadays the blow-up method, is a
generalization of the classical Harnack construction of $M$-curves
by small perturbations, using the components of the exceptional
divisor as a basis of rank one (see \cite{RisH}, \cite{KR} and
\cite{KRS}). One of our motivations was to study to which extent
Mikhalkin's result holds for smoothings of singular points of real
plane curves, particularly for {\em Harnack smoothings}, those
$M$-smoothings which are in maximal position with respect to two
coordinates lines through the singular point.

We develop a new construction of smoothings of a real plane branch
$(C,0)$ by using Viro Patchworking method. Since real plane
branches are Newton degenerated, we cannot apply Viro Patchworking
method directly. Instead we apply the Patchworking method for
certain {\em Newton non degenerate} curve singularities with
several branches which are defined by semi-quasi-homogeneous
polynomials. These singularities appear as a result of iterating
deformations of the strict transforms $(C^{(j)}, o_j)$, at certain
infinitely near points $o_j$ of the embedded resolution of
singularities of $(C,0)$. Our method yields multi-parametric
deformations, which we call {\em multi-semi-quasi-homogeneous}
(msqh) and provides simultaneously msqh-smoothings of  the strict
transforms $(C^{(j)}, o_j)$. We exhibit suitable hypothesis which
characterize $M$-smoothings and Harnack smoothings for this class
of deformations (see Theorem \ref{main}). Up to the author's
knowledge, Theorem \ref{main} is the first instance in the
literature in which Viro Patchoworking method is used to define
smoothings of Newton degenerated singularities, with controlled
topology.

We introduce the notion of {\em multi-Harnack smoothings}, those
Harnack smoothings, such that the msqh-$M$-smoothings of the
strict transforms  $C^{(j)}$ appearing in the process are Harnack.
We prove that any real plane branch $C$ admits a multi-Harnack
smoothing. For this purpose we prove the existence of Harnack
smoothings of singularities defined by certain
semi-quasi-homogeneous polynomials (see Proposition
\ref{compatible}). One of our main results, Theorem
\ref{main-top}, states that multi-Harnack smoothings of a real
plane branch $(C,0)$ have a unique topological type, which depends
only on the complex equisingular class of $(C,0)$.  In particular,
multi-Harnack smoothings do not have nested ovals. Theorem
\ref{main-top} can be understood as a local version of Mikhalkin's
Theorem \ref{Mikhalkin}. The proof is based on Theorem \ref{main}
and an extension of Mikhalkin's result for Harnack smoothings of
certain non-degenerated singular points (Theorem \ref{1par}).

We also analyze certain multi-scaled regions containing the ovals.
The phenomena is quite analog to the analysis of the asymptotic
concentration of the curvature of the Milnor fibers in the complex
case, due to Garc\' \i a Barroso and Teissier \cite{GBT}.

It is a challenge for the future to extend, as possible, the
techniques and results of this paper to the constructions of
smoothings of other singular points of real plane curves.

The paper is organized as follows. The five first sections are
preliminary material: in Sections \ref{definitions}, \ref{toric}
and \ref{ViroT} we introduce the Viro Patchworking method, also in
the toric context;  we recall Mikhalkin's result on Harnack curves
in projective toric surfaces in Section \ref{Mik}; the notion of
smoothings of real plane curve singular points is presented in
section \ref{smoothings}. Section \ref{harnack-smoothings}
contains the first new results, in particular, the determination
of the topological type of Harnack smoothings of singularities
defined by certain non degenerated semi-quasi-homogeneous
polynomials (see Theorem \ref{1par}). In Section \ref{toric-res}
we recall the construction of a toric resolution of a plane branch
and we introduce the support of the msqh-smoothings, which are
studied in the last section. The main results of the paper are
collected in Section \ref{last}: the characterization of maximal,
Harnack and multi-Harnack-msqh-smoothings in Theorem \ref{main}
and Corollary \ref{multi-h2}, and the characterization of the
topological type of multi-Harnack smoothings in Theorem
\ref{main-top}, the description of the scales of   ovals in
Section \ref{taille}  and finally some examples explained in
detail.

\section{Basic notations and definitions} \label{definitions}

A real algebraic variety is a complex algebraic variety $V$
invariant by an anti-holomorphic involution, we denote by $\R V$
the real part. For instance, a real algebraic plane curve $C
\subset \C^2$ is a complex plane curve which is invariant under
complex conjugation.  The curve $C$ is defined by $P =0$ where $0
\ne P = \sum c_{i,j}  x^i y^j \in \R [x, y]$. We use the following
notations and definitions:

The {\em Newton polygon} of $P$  (resp. the {\em local Newton
polygon}) is the convex hull in $\R^2$ of the set $\{ (i, j) \mid
c_{i, j} \ne 0 \}$ (resp. of  $\{ (i, j) + \R^2_{\geq 0} \mid
c_{i,j} \ne 0 \}$).

 If $\Lambda \subset \R^2$ we
denote by $P^{\Lambda} $ the symbolic restriction $P^{\Lambda}:=
\sum_{(i, j) \in \Lambda \cap \Z^2} c_{i,j} x^i y^j$.

Suppose that  $0 \in C \subset \C^2$ is an isolated singular point
of $C$ and that $C$ does not contain a coordinate axis. Then the
Newton diagram of $P$ is the closed region $\D$ bounded by the
coordinate axis and the local Newton polygon of $P$.

The polynomial $P$ is {\em non degenerated} (resp. {\em real non
degenerated}) with respect to its Newton polygon if for any
compact face $\Theta$ of it we have that $P^{\Theta} =0 $ defines
a non singular subset of $(\C^*)^2$ (resp. of $(\R^*)^2$). In this
case if $\Lambda$ is an edge of the Newton polygon of $P$, then
the polynomial $P^\Lambda$ is of the form:
\begin{equation} \label{peripheral}
P^{\Lambda} = c x^a y ^b \prod_{i=1}^e (y^{n} - \a_i x^{m}),
\end{equation}
where $c \ne 0$, $a, b \in \Z_{\geq 0} $, the integers $n, m \geq
0$ are coprime and the numbers $\a_1, \dots,\a_e \in \C^*$, which
are called {\em peripheral roots of $P$ along the edge $\Lambda$
or simply peripheral roots of $P^{\Lambda}$}, are distinct (resp.
the real peripheral roots of $P^{\Lambda}$ are distinct). Notice
that the non degeneracy (resp. real non degeneracy) of $P$ implies
that $C \cap (\C^*)^2$ is non singular (resp. $C \cap (\R^*)^2$),
for taking $\Theta$ equal to the Newton polygon of $P$ in the
definition.

We say that $P$ is {\em non degenerated} (resp. {\em real non
degenerated}) with respect to its local Newton polygon if for any
edge $\Lambda$ of it we have that  the equation $P^{\Lambda} =0 $
defines a non singular subset of $(\C^*)^2$ (resp. of $(\R^*)^2$).

The notion of non degeneracy with respect to the Newton polyhedra
extend for polynomials of more than two variables (see
\cite{Kou}).

\section{The real part of a projective toric variety} \label{toric}

We introduce basic notations and facts on the geometry of toric
varieties. We refer the reader to \cite{GKZ}, \cite{Oda} and
\cite{Fulton} for proofs and more general statements. For
simplicity we state the notations only for surfaces.

Let $\Theta$ be a convex two dimensional polytope in $\R^2_{\geq
0}$ with integral vertices, a {\em polygon} in what follows. We
associate to the polygon $\Theta$ a projective toric variety
$Z(\Theta)$. The algebraic torus $(\C^*)^2$ is embedded as an open
set of $Z(\Theta)$, and acts on $Z(\Theta)$, in a way which
coincides with the group operation on the torus. There is a one to
one correspondence between the faces of $\Theta$ and the orbits of
the torus action, which preserves the dimensions and the
inclusions of the closures. If $\Lambda$ is a one dimensional face
of $\Theta$, we have an embedding $Z(\Lambda) \subset Z(\Theta)$.
The variety $Z(\Lambda)$ is a projective line $\C P^1$ embedded in
$Z(\Theta)$. These lines are called the {\em coordinate lines} of
$Z(\Theta)$. The intersection of two coordinate lines
$Z(\Lambda_1)$ and $Z(\Lambda_2)$ reduces to a point (resp. is
empty) if and only if the edges $\Lambda_1$ and $\Lambda_2$
intersect in a vertex of $\Theta$ (resp. otherwise). The surface
$Z(\Theta)$ may have singular points only at the zero-dimensional
orbits. The algebraic real torus $(\R^*)^2 $ is an open subset of
the real part $\R Z( \Theta)$  of $Z (\Theta)$ and acts on it. The
orbits of this action are just the real parts of the orbits for
the complex algebraic torus action. For instance if $\Theta$ the
simplex with vertices  $(0,0)$, $(0,d)$ and $(d,0)$ then the
surface $Z (\Theta)$ with its coordinate lines
 is the complex projective plane with the classical three
coordinate axis.

The  image of the {\em moment map} $\f_{\Theta} : (\C^*)^2
\rightarrow \R^2$,
\begin{equation}\label{moment}
(x, y) \mapsto \left(\sum_{\a_k ={(i,j)}, \a_k \in \Theta \cap
\Z^2} |x^{i} y ^{j}| \right)^{-1} \left( \sum_{\a_k ={(i,j)}, \a_k
\in \Theta \cap \Z^2} | x^{i} y ^{j}| (i, j) \right).
\end{equation}
is $\mbox{{\rm int}} \Theta$, where $\mbox{{\rm int}}$ denotes
relative interior. The restriction $\f_{\Theta}^+ := \f_{\Theta |
\R^2_{>0}} $ is a diffeomorphism of $ \R^2_{>0}$ onto the interior
of $\Theta$.

Denote by ${\mathcal S} \cong (\Z/ 2 \Z)^2$ the group consisting
of the orthogonal  symmetries of $\R^2$ with respect to the
coordinate lines, namely the elements of ${\mathcal S}$ are
$\r_{i,j} : \R^2 \rightarrow \R^2 $, where $\r_{i,j} (x, y) =
((-1)^i x, (-1)^j y)$ for $(i, j) \in \Z^2_2$.

If $A \subset \R^2$ we denote by $\tilde{A}$ the union $\tilde{A}
:= {\bigcup_{\r \in {\mathcal S}}}  \,  \r( A) \subset \R^2$.

The map $\f_{\Theta}^+$  extends to a diffeomorphism: $
\tilde{\f}_{\Theta} : (\R^*)^2 \rightarrow \widetilde{\mbox{{\rm
 int}}{({\Theta})}}$ by $\tilde{\f}_{\Theta} ( \r (x) ) := \r
( \f( x))$, for  $x \in \R^2_{>0}$ and $\r \in {\mathcal S}$.

If $\Lambda$ is an edge of $\Theta$ and if $n= (u, v) $ is a
primitive integral vector orthogonal to $\Lambda$ we denote by
$\r_\Lambda$ or by $\r_\Lambda (a,b)$ the element of ${\mathcal
S}$ defined by $\r_\Lambda (a,b) = ( (-1)^{u} a, (-1)^{v} b)$.

We consider the equivalence relation $\sim$ in the set
$\tilde{\Theta}$, which for each edge $\Lambda$ of $\Theta$,
identifies a point in $\Lambda$ with its symmetric image by
$\r_\Lambda$. Set $\tilde{\Lambda}/\sim$ the image of
$\tilde{\Lambda}$ in $\tilde{\Theta}/\sim$. For each edge
$\Lambda$ of $\Theta$ we have diffeomorphisms ${\f}_{\Lambda}^+ :
\R_{>0}  \rightarrow \mbox{\rm int} {\Lambda}$ and
$\tilde{\f}_{\Lambda} : \R^*  \rightarrow \widetilde{\mbox{\rm
int} {\Lambda}}/\sim $, corresponding to the moment map in the one
dimensional case. Notice that $\widetilde{\mbox{\rm int}
{\Lambda}}/\sim $ has two connected components  and the real part
$\R Z(\Lambda)$ corresponds to $\tilde{\Lambda}/\sim$.

We summarize this constructions in the following result (see
\cite{Oda} Proposition 1.8 and \cite{GKZ}, Chapter 11, Theorem 5.4
for more details and precise statements).

\begin{Pro} \label{quatro}
The morphisms defined from the moment map glue up ia a stratified
homeomorphism
\begin{equation} \label{psi-theta}
\begin{array}{cccc}
\Psi_{\Theta} :  & \tilde{\Theta}/\sim \,  &  \longrightarrow &
\R Z(\Theta),
\end{array}
\end{equation}
which for any edge $\Lambda$ of $\Theta$ applies
$\tilde{\Gamma}/\sim$ to the correspondent coordinate line $\R
Z(\Lambda) \subset \R Z( \Theta)$. The composite
\[
\begin{array}{ccccccc}
(\R^*)^2  & \stackrel{\tilde{\f}_\Theta}{\longrightarrow} &
\tilde{\Theta} &  \longrightarrow &   \tilde{\Theta} /\sim
 & \stackrel{ \Psi_{\Theta} }{\longrightarrow} & \R Z (\Theta).
\end{array}
\]
is the inclusion of the real part of the torus in $\R Z (\Theta)$.
\end{Pro}


 A polynomial $P
\in \R [x, y]$ with Newton polygon equal to $\Theta$ defines a
real algebraic curve $C$ in the real toric surface $Z(\Theta)$
which does not pass through any $0$-dimensional orbit. If $C$ is
smooth its genus coincides with the number of integral points in
the interior of $\Theta$, see \cite{Kh}. The curve $C$ is a
$M$-curve if $\R C$ has the maximal number $ 1 + \# \mbox{{\rm
int}} {\Theta} \cap \Z^2$ of connected components. If $\Lambda$ is
an edge of $\Theta$ the intersection of $C$ with the coordinate
line $Z(\Lambda)$ is defined by $P^{\Lambda}$.  The number of
zeroes of $P^{\Lambda}$ in the projective line $Z(\Lambda)$,
counted with multiplicity, is equal to the {\em integral length}
of the segment $\Lambda$, i.e., one plus the number of integral
points in the interior of $\Lambda$. This holds since these zeroes
are the image of the peripheral roots $\a_i$ of $P^\Lambda$ (see
(\ref{peripheral})) by the embedding map $\C^* \hookrightarrow
Z(\Lambda)$.
 For
this reason we abuse sometimes of terminology and call peripheral
roots the zeroes of $P^{\Lambda}$ in the projective line
$Z(\Lambda)$.

\section{Patchworking real algebraic curves} \label{ViroT}

Patchworking  is a method introduced by Viro for constructing real
algebraic hypersurfaces  (see \cite{Vth}, \cite{Vpw} \cite{V67},
\cite{VirL} and  \cite{IV}, see also \cite{GKZ} and \cite{RisV}
for an exposition and \cite{Bihan} and \cite{St} for some
generalizations). We use the Notations introduced in Section
\ref{toric}.

Let $\Theta \subset \R^2_{\geq 0}$ be an integral polygon. The
following  definition is fundamental.

\begin{Defi} \label{chart}
Let $Q \in \R[x, y]$ define a real algebraic curve $C$ then  the
$\Theta$-chart $\mbox{\rm Ch}_\Theta (C)$ of $C$ is the closure of
$\mbox{\rm Ch}_\Theta^* (C) := \tilde{\f}_\Theta ( \R C  \cap
(\R^*)^2 ) $ in $\tilde{\Theta}$.
\end{Defi}

 If $\Theta $  is the Newton polygon of $Q$ we often  denote
$\mbox{\rm Ch}_\Theta (C)$ by $\mbox{\rm Ch}(Q)$ or by $\mbox{\rm Ch}(C)$
if the coordinates used are clear from the context. If $Q$ is real
non degenerated with respect its Newton polygon, then for any face
 $\Lambda$ of
 $\Theta$ we have that $\mbox{\rm Ch}_\Theta (Q) \cap
\tilde{\Lambda} = \mbox{\rm Ch}_\Lambda ( Q_\Lambda)$ and the
intersection $\mbox{\rm Ch}_\Theta (Q) \cap \tilde{\Lambda}$ is
transversal (as stratified sets).

\begin{Not} \label{notation-viro}
We consider a polynomial $P_t = \sum A_{i,j}  (t)  x^i y^j \in \R
[t, x, y]$, as a family of polynomials in $\R[x, y]$.
\begin{enumerate}
\item We denote by $\Theta \subset \R^2$ the Newton polygon of
$P_t$, when $0< t \ll 1$. \item We denote by $\hat{\Theta} \subset
\R \times \R^2$ the Newton polytope of $P_t$, when it is  viewed
as a polynomial in $\R[t, x, y]$.
\item We denote by $
\hat{\Theta}_{c} $ the {\em lower part} of $\hat{\Theta}$, i.e.,
the union of compact faces of the Minkowski sum $\hat{\Theta} +
(\R_{\geq 0} \times \{ (0,0) \}) $.

\item The restriction of the second projection $\R \times \R^2
\rightarrow \R^2  $ to $\hat{\Theta}_{c} $ induces a finite
strictly convex polyhedral subdivision $\Theta'$ of $\Theta$. The
inverse function $ \w: \Theta \rightarrow \hat{\Theta}_{c} $ is a
{\em piece-wise affine strictly convex function}. Any cell
$\Lambda$ of the subdivision $\Theta'$ corresponds by this
function to a face $\hat{\Lambda}$ of $\hat{\Theta}$  contained in
$\hat{\Theta}_{c}$,  of the same dimension,  and the converse also
holds. The Newton polygon of $P_0$ is a face of $\Theta'$ by
construction.
 \item If $\Lambda$ is a cell of $\Theta'$ we denote
by $P^{\hat{\Lambda}}_1$, or by $P^{\hat{\Lambda}}_{t=1}$ the
polynomial in $\R[x,y]$ obtained by substituting $t=1$ in
$P_t^{\hat{\Lambda}}$.
\end{enumerate}
\end{Not}

\begin{The} \label{viro} With the above notations,
if for each face  $\hat{\Lambda}$ of $\hat{\Theta}$ contained in
$\hat{\Theta}_{c} $ the polynomial $P^{\hat{\Lambda}}_{{1}}$ is
real non degenerated with respect to ${\Lambda}$, then the pair $(
\tilde{\Theta},  \mbox{\rm Ch}_\Theta (P_t)  )$ is stratified
homeomorphic to $(\tilde{\Theta}, \tilde{C})$, where $\tilde{C}$
is the curve obtained by gluing together in $\tilde{\Theta}$ the
charts $\mbox{\rm Ch}_\Lambda ( P^{\hat{\Lambda}}_{{1}} ) $ for
$\Lambda$ running through the cells of the subdivision $\Theta'$,
for $0 <t \ll 1$.
\end{The}

\begin{Rem} \label{v-gen}
The statement of Theorem \ref{viro} above is a slight
generalization of the original result of Viro in which the
deformation is of the form $ \sum  A_{i,j} t^{\w(i,j)} x^i y^j $,
for some real coefficients $A_{i,j}$. The same proof generalize
without relevant changes to the case presented here, when we may
add terms $B_{k, i,j} t^k x^i y^j$ with $B_{k, i,j} \in \R$ and
exponents $(k, i, j)$ contained in $(\hat{\Theta} + \R_{\geq 0}
\times \{ 0\} ) \setminus \hat{\Theta}_{c} $.
\end{Rem}

\begin{Rem} \label{viro-toric}
Theorem \ref{viro} extends naturally to provide constructions of
real algebraic curves with prescribed topology in the real toric
surface $Z (\Theta)$. The chart $\mbox{Ch}_{\tilde{\Theta} /\sim}
(C) $ of the curve $C$ in $Z (\Theta) $ is defined as the closure
of $\mbox{Ch}^*_{\Theta} (C) $ in $\tilde{\Theta}/\sim$ (where
$\sim$ is the equivalence relation defined in Section
\ref{toric}). Then the statement of Theorem \ref{viro} holds for
the curve $C_t$ defined by $P_t$ in $Z(\Theta)$ by identifying
$\tilde{\Theta}/\sim$ and $\R Z( \Theta)$ by the map $\Psi_\Theta$
(see (\ref{psi-theta})).
\end{Rem}

\begin{Defi}
The gluing of charts of Theorem \ref{viro} is called {\em
combinatorial patchworking} if the subdivision $\Theta'$ is a {\em
primitive triangulation} of $\Theta$, i.e., it contains all
integral points of $\Theta$ as vertices.
\end{Defi}
Notice that $\Theta'$ is a primitive triangulation if  the two
dimensional cells $\Lambda$ are primitive triangles, i.e., of area
$1/2$ with respect to the standard volume form induced by a basis
of the lattice $\Z^2$. The description of the charts of a
combinatorial patchworking is determined by the subdivision
$\Theta'$ and the signs of the terms appearing in $P_{t}^{
\hat{\Theta}_c} $ as a polynomial in $x$ and $y$. The  {\em
distribution of signs}  $\epsilon: \Theta \cap \Z^2 \rightarrow \{
\pm 1\}$,  induced by taking the signs of the terms appearing in
$P_{t}^{ \hat{\Theta}_c} $, extends to
\begin{equation} \label{tilde}
\tilde{\epsilon} : \tilde{\Theta} \cap \Z^2 \rightarrow \{ \pm
1\}, \mbox{ by  setting } \tilde{\epsilon} (r,s) = (-1)^{i+j}
\epsilon \circ \r_{i,j} (r,s) \mbox{ whenever } \r_{i,j} (r,s) \in
\Z^2_{\geq 0}.
\end{equation}
The chart associated to the polynomial $P_t$ in a triangle is
empty if all the signs are equal and otherwise is isotopic to
segment dividing the triangle in two parts, each one containing
only vertices of the same sign. See \cite{GKZ}, \cite{IV},
\cite{Itenberg-T}.

\section{Maximal and Harnack curves in projective toric
surfaces} \label{Mik}

If $C $ is a smooth real projective curve of degree $d$ then the
classical Harnack inequality states that the number of connected
components of its real part $\R C \subset \R P^2 $ is bounded by
$\frac{1}{2} (d-1)(d-2) +1$. The curve $C$ is called maximal or
$M$-curve if the number of components is equal to the bound.
Maximal curves always exists and geometric constructions of such
curves were found in particular by Harnack, Hilbert and Brusotti.
Determining the possible topological types of the pairs $(\R C, \R
P^2)$ in terms of the degree $d$ is usually called the first part
of the Hilbert's 16th problem.
\begin{Defi} \label{max-pos}
A real projective curve $C$ of degree d is in
\begin{enumerate}
\item[(i)] {\em maximal position} with respect to a  real line $L$
if the intersection $L\cap C$ is transversal, $L\cap C = \R L \cap
\R C$  and  $L \cap C$ is contained in one connected component of
$\R C$.

\item[(ii)] {\em maximal position with respect to real lines}
$L_1, \dots, L_n$  in $\C P^2$ if $C$ is in maximal position with
respect to $L_i$, and  there exist $n$ disjoint arcs
$\mathbf{a}_1, \dots, \mathbf{a}_n$ contained in one connected
component of $\R C$ such that $C \cap L_i = \mathbf{a}_i \cap \R
L_i$, for $i=1, \dots, n$.
\end{enumerate}
\end{Defi}
  Mikhalkin studied the topological types of the triples
$(\R P^2, \R C, \R L_1\cup \dots \cup \R L_{n})$ for those
$M$-curves $C$ in maximal position with respect to lines $L_1,
\dots, L_n$. He proved that for $n = 3$ there is a unique
topological type, while for $n> 3$ there is none. For $n=1$ the
classification reduces to the topological classification of
maximal curves in the real affine plane (which is open for $d >
5$), while for $n=2$ there are several constructions of $M$-curves
of degree $d$, which were found by Brusotti, with $d \geq 4$. See
\cite{Brusotti} and \cite{Mikhalkin}.

Mikhalkin's results were stated and proved  more generally for
real algebraic curves in projective toric surfaces. We denote by
$L_1, \dots, L_m$ the sequence of cyclically incident coordinate
lines in the toric surface $Z(\Theta)$ associated to the polygon
$\Theta$. The notion of maximal position of  a real algebraic
curve $C \subset Z(\Theta)$ with respect one line is the same as
in the projective case.

\begin{Defi}
A real algebraic curve $C$ in the toric surface $Z(\Theta)$ is in
{\em maximal position} with respect to lines
 $L_1, \dots, L_n$, for $1
\leq n \leq m$, if there exist $n$ disjoint arcs $\mathbf{a}_1,
\dots, \mathbf{a}_n$ contained in one connected component of $\R
C$
 such that the
intersection   $C \cap L_i$ is transversal and contained in
$\mathbf{a}_i \cap \R L_i$, while $\mathbf{a}_i \cap \R L_j =
\emptyset$ if $i \ne j$, for $ i =1, \dots, n$ and $j =1, \dots,
m$. In addition, we say that:
\begin{enumerate}
\item[(i)] The curve $C$ is {\em cyclically in maximal position}
if it is in maximal position with respect to the lines $L_1,
\dots, L_m$ and the points of intersection of $\R C$ with the
lines $L_1, \dots, L_m$  when viewed in the connected component of
$\R C$ are partially ordered, following the adjacency of the lines
$L_1, \dots, L_m$ .
 \item[(ii)] The curve
 $C$ has {\em good oscillation} with respect to the line $L_i$ if the
points of intersection of $C$ with $L_i$ have the same order when
viewed in the arc $\mathbf{a}_i$ and in the line $\R L_i$.
\end{enumerate}
\end{Defi}

\begin{Rem}
If $C \subset Z (\Theta)$ is in maximal position with respect to
the coordinate line $L = Z(\Gamma)$, for $\Gamma$ an edge of
$\Theta$, we say also that the chart $\mbox{\rm Ch}_{\Theta/\sim}
(C)$ is in maximal position with respect to $\tilde{\Gamma}/\sim$
(see Proposition \ref{quatro}).
\end{Rem}

We have the following result of Mikhalkin (see \cite{Mikhalkin}).
\begin{The} \label{Mikhalkin}
With the above notations if a $M$-curve $C$ is cyclically in
maximal position with respect to the coordinate lines $L_1, \dots,
L_m$ of the real toric surface $ Z( \Theta)$ then the topological
type of the triple $(\R C, \R Z (\Theta), (\R^*)^2)$ depends only
on $\Theta$.
\end{The}

\begin{Defi}
A Harnack curve in the real toric surface $Z(\Theta)$ is a real
algebraic curve  $C$ verifying the conclusion of Theorem
\ref{Mikhalkin}.
\end{Defi}

\begin{Rem}\begin{enumerate}
\item The notion of Harnack curve in this case depends on the
polygon $\Theta$. Changing $\Theta$ by $k \Theta$, for $k >1$
provides the same toric variety but the corresponding Harnack
curves are different.

\item By our convention and the definitions, the Newton polygon of
a polynomial $P \in \R[x, y]$, defining a Harnack curve $C \subset
Z (\Theta)$ is equal to $\Theta$. Notice also that $P$ is non
degenerated with respect to $\Theta$ and for any edge $\Lambda$ of
$\Theta$ the peripheral roots of $P$ along the edge $\Lambda$ are
real and of the same sign (see (\ref{peripheral})).
\end{enumerate}
\end{Rem}

 If $C$
defines a Harnack curve in $Z( \Theta)$ we denote by $\O_C$ the
unique connected component of $\R C$ which intersects the
coordinate lines and by
 $U_C$ the set $\cup {B}_k$, where $B_k$ runs through the connected components of the set
 $(\R^*)^2 \setminus \O_C$, whose boundary meets at most two
 coordinate lines. The set $U_C$ is an open region
 bounding $\O_C$ and the coordinate lines of the toric surface.

 The following Proposition, which is a reformulation of the
 results of Mikhalkin \cite{Mikhalkin}, describes completely the
 topological type of the real part of a Harnack curve in a toric
 surface. See Figure \ref{normalizado} below an example.
\begin{Pro} \label{top-harnack}
With the previous notations,  if $C$ defines a Harnack curve in
$Z( \Theta)$ then $C$ is cyclically in maximal position with good
oscillation with respect to the coordinate lines of the toric
surface $Z(\Theta)$. Moreover, the set of connected components of
 $\R C \cap (\R^*)^2$ can be labelled as $\{ \O_{r,s} \}_{(r,s)
 \in \Theta \cap \Z^2}$ in such a way that:

 - there exists a unique $(i, j) \in \Z^2_2$
such that for any $(r,s) \in \Theta \cap Z^2$:
\begin{equation} \label{indice}
\O_{r,s} \subset  \r_{i,j} ( \R^2_{s,r} ).
\end{equation}

- the set of components $\{ \O_{r,s} \}_{(r,s) \in \mbox{{\rm
\small int}} \Theta}$ consists of non nested ovals in $(\R^*)^2
\setminus U_C$.

- if $(r,s), (r',s') \in \partial \Theta \cap \Z^2$ then the
intersection of closures $\overline{\O}_{(r,s)} \cap
\overline{\O}_{(r',s')}$  reduces to a point in $Z(\Lambda)$ if
and only if  $(r,s)$ and $(r',s')$ are consecutive integral points
in $\Lambda$ for some edge $\Lambda$ of $\Theta$, or it is empty
otherwise.

- we have that  $\O_C = \bigcup_{(r,s) \in
\partial \Theta \cap \Z^2} \overline{\O}_{r,s}$.
\end{Pro}

The following Proposition describes the construction of  Harnack
curve in the real toric surface $Z(\Theta)$ by using combinatorial
patchworking (see \cite{Itenberg-T} Proposition 3.1 and
\cite{Mikhalkin} Corollary A4).

\begin{Defi} \label{harnack-sign}
Denote by $\epsilon_h : \Z^2 \rightarrow  \{ \pm 1\}$ the function
$\epsilon_h (r,s) = (-1)^{(r-1)(s-1)}$. A {\em Harnack
distribution of signs} is any of the following distributions $ \pm
\tilde{\epsilon}_h \circ \r_{i, j} $, where  $(i, j) \in \{0,
1\}^2$ (see (\ref{tilde})).
\end{Defi}

\begin{Pro} \label{harnack-curve}
Let $\epsilon$ be a Harnack distribution of signs and $\w: \Theta
\cap \Z^2 \rightarrow \Z$ define a primitive triangulation of
$\Theta$, then the polynomial
\[
P_t = \sum_{(i,j) \in \Theta \cap \Z^2} \epsilon (i, j)
t^{\w(i,j)} x^i y^j
\]
defines a Harnack curve in $Z(\Theta)$, for $0 < t \ll 1$.
\end{Pro}
{\em Proof.} Suppose without loss of generality that $\epsilon =
\epsilon_h$. Since the triangulation is primitive any triangle $T$
containing a vertex with both even coordinates, which we call even
vertex, has the other two vertices with at least one odd
coordinate. It follows that the even vertex has a sign different
than the two other vertices. If the even vertex is in the interior
of $\Theta$ then there is necessarily an oval around it, resulting
of the combinatorial patchworking. The situations in the other
quadrants is analogous by the symmetry of the Harnack distribution
of signs since, for each triangle $T$, for any vertex $v \in T$
there is a unique symmetry $\r_{i,j}$ such that the sign of
$\r_{i,j} (v)$ is different than the sign of the two other
vertices of $\r_{i,j} (T)$. It follows that there are $\#
\mbox{{\rm int}} \Theta \cap \Z^2$ ovals which do not cut the
coordinate lines and exactly one more component which has good
oscillation with maximal position with respect to the coordinate
lines of $Z(\Theta)$. $\Box$

\section{Smoothings of real plane singular points}
\label{smoothings}

Let $(C,0)$ be a germ of real plane curve with an isolated
singular point at the origin.  Set $B_{\epsilon} (0)$  for the
ball
 of center $0$ and of radius $\epsilon$. If $0 < \epsilon \ll 1$
each branch of $(C,0)$ intersects $\partial B_{\epsilon} (0)$
transversally along a smooth circle and the same property holds
when the radius is decreased (analogous statements hold also for
the real part). Then we denote the ball $B_\epsilon (0)$ by
$B(C,0)$, and we called it a {\em Milnor ball} for $(C,0)$. See
\cite{Milnor}).

A smoothing $C_t  \subset B_\epsilon (C, 0)$ is a real analytic
family $C_t \subset B(C,0)$, for $t \in [0,1]$ such that $C_0 = C$
and $C_t$ for $0 < t \ll 1 $ is non singular and transversal to
the boundary. By the connected components of the smoothing $C_t$
we mean the components of the real part $\R C_t$ of $C_t$ in the
Milnor ball. These components consist of finitely many ovals and
non closed connected components.

\begin{Defi}   \label{topological2}
\begin{enumerate}\item[(i)] The {\em topological type}  of a smoothing $C_t$ of a
plane curve singularity $(C, 0)$ with Milnor ball $B= B(C,0)$ is
the topological type of the pair $(\R C_t \cap \R {B}, \R {B})$.

\item[(ii)] The {\em signed topological type}  of a smoothing
$C_t$ of a plane curve singularity $(C, 0)$ with Milnor ball $B=
B(C,0)$ with respect to fixed coordinates $(x, y)$ is the
topological type of the pairs,  $(\R C_t \cap \R {B}, \R {B} \cap
\R^2_{i, j})$, for $(i, j) \in \{0, 1 \}^2$, where $\R^2_{i, j} $
denotes the open quadrant $\R^2_{i,j} := \{ (x, y ) \mid (-1)^i x
> 0, (-1)^j y > 0 \}$.
\end{enumerate}
\end{Defi}

A branch of $(C,0)$ is {\em real} if it has a Newton Puiseux
parametrization with real coefficients. Denote by $r$ the number
of (complex analytic) branches of $(C,0)$. The number of non
closed components is equal to $r_\R$, the number of {\em real
branches} of $(C, 0)$. The number of ovals of a smoothing $C_t$ is
$\leq \frac{1}{2} ( \mu (C)_0 - r +1 )$ if $r_\R \geq 1$ and $\leq
\frac{1}{2} ( \mu (C)_0 - r + 3)$ if $r_\R =0$; where $\mu (C)_0$
denotes the Milnor number of $C$ at the origin (see \cite{Arnold},
\cite{RisH}, \cite{KOS} and \cite{KR}). A smoothing is called a
$M$-smoothing if the number of ovals is equal to the  bound.

The existence of $M$-smoothings is a quite subtle problem, for
instance if $(C,0)$ is a real plane branch there exists always a
$M$-smoothing (see \cite{RisH}) however there exists singularities
which do not have a $M$-smoothing (see \cite{KOS}). Some other
types of real plane singularities which do have a $M$-smoothing
are described in \cite{KR} and \cite{KRS}.

\subsection{Patchworking smoothings  of plane curve singularities}
\label{smoothing}

With the notations as above we consider a family of polynomials
$P_t (x,y) $  such that  $P_0(0,0) =0$, $\mbox{\rm ord}_x
P_0(x,0)= m>1$, $\mbox{\rm ord}_y P_0 (0, y) = n >1$ and $P_t (0,
0) \ne 0$ defining a smoothing $C_t$ of the germ of plane curve
singularity $(C,0)$ of equation $P_0 (x,y) =0$. Then the germ
$(C,0)$ does not contain any of the coordinate axis and  $0 \notin
C_t$ for $t \ne 0$. The Newton diagram $\D$ of $P_0$ is contained
in the Newton polygon $\Theta$ of $P_t$, for $0 < t \ll 1$.

The following result is a consequence of Theorem \ref{viro}, see
\cite{VirL},  \cite{Vpw} and \cite{KRS}.

\begin{The} \label{vlocal} With Notations \ref{notation-viro},
the family $P_t$ defines a
convex subdivision $\Theta'$ of $\Theta$ in polygons. If
$\Lambda_1, \dots, \Lambda_k$ are the cells of $\Theta'$ contained
in $\D$ and if
 $P^{\hat{{\Lambda}}_i}_{{1}}$ is real non degenerated
with respect to ${\Lambda}_i$ for $i=1, \dots, k$, then the family
$P_t$ defines a smoothing $C_t$ of $(C,0)$ such that the pairs
$({\R {B} (C,0)}, \R C_t)$ and $(\tilde{\D}, \bigcup_{i=1}^k
\mbox{\rm Ch}_{\Lambda_i} (P^{\hat{{\Lambda}}_i}_{{1}}))$ are
homeomorphic (in a stratified sense), for $0 <t \ll 1$.
\end{The}

The hypothesis of Theorem \ref{vlocal} imply that $P_0$ is real
non degenerated with respect to its local Newton polygon.

\subsection{Semi-quasi-homogeneous smoothings} \label{Section-sqh}

We say that a polynomial  $P_0 \in \R [x,y]$, with $P(0) =0$  is
 {\em semi-quasi-homogeneous} (sqh)
if its local Newton polygon has only one compact edge.

\begin{Not} \label{notation-sqh}
Consider a sqh-polynomial $P_0 \in \R[x, y]$ with Newton polygon
with compact edge $\Gamma = [(m,0), (0, n)]$  and Newton diagram
$\D=[(m,0), (0, n), (0,0)]$. We denote by $\D^-$ the set $\D
\setminus \Gamma$. The number $e := \mbox{\rm gcd} (n, m) \geq 1$
is equal to the
 {integral length} of the segment $\Gamma$. We set $n_0 = n/e $
and $m_0 = m/e$. Notice that the polynomial $P_0$ is of the
following form:
\begin{equation} \label{edge-2}
P_0  = \prod_{s =1}^e (y^{n_0} - \vartheta_{s} x^{m_0}) + \cdots,
\mbox{ for some } \vartheta_{s} \in \C^*,
\end{equation}
 where the exponents $(i, j)$ of the terms which are not
written verify that $n i + mj > nm$.
 Denote
by $\D^-$ the set $\D \setminus \Gamma$.
\end{Not}

Suppose that all the peripheral roots $\vartheta_s$  of $P^\Gamma$
are real and different. In this case, using Kouchnirenko's
expression for the Milnor number of $(C, 0)$ (see \cite{Kou}), we
deduce that the bound on the number of connected components (resp.
ovals) of a smoothing of $(C,0)$ is equal to
\begin{equation} \label{mu-delta}
\# ( \mbox{\rm int}  {\D} \cap \Z^2) + e -1, \quad ( \mbox{\rm
resp.} \# ( \mbox{\rm int}  {\D} \cap \Z^2) \, ).
\end{equation}

We consider a uniparametrical family of polynomials $P_t (x, y)$
with real coefficients with $P_t (0,0) \ne 0$ and $P_0 (x, y)$ of
the form (\ref{edge-2}).

Consider the lower part $\hat{\Theta}_c$ of the Newton polyhedra
of $P_t$ viewed as a polynomial in $\R[t, x, y]$. The projection
of the faces of $\hat{\Theta}_c$ define a polygonal subdivision
$\Theta'$ of the Newton polygon $\Theta \subset \R^2$ of $P_t$,
viewed as a polynomial in the variables $x$ and $y$.  See
Notations \ref{notation-viro}.

\begin{Defi} \label{sqh}
We say that $C_t$  is a semi-quasi-homogeneous (sqh) deformation
(resp. smoothing)  of $(C,0)$ if $\D$ is a face of the subdivision
$\Theta'$ (resp. and in addition the polynomial
$P_{{1}}^{\hat{\D}}$ is real non degenerated).
\end{Defi}
Notice that if $C_t$ is a sqh-deformation the polynomial
$P_t^{\hat{\D}}$ is quasi-homogeneous as a polynomial in $(t,
x,y)$.  This implies that any non zero monomial of
$P_{{t}}^{\hat{\D}}$ is of the form $a_{i,j} t^{w_{i,j}} x^i y^j$
where the ratio of $w_{i,j}$ by $nm  - n i - m_j$ is some positive
constant for all $i, j \in \D^-$.

If $P_{{1}}^{\hat{\D}}$ is real non degenerated then the real
peripheral roots of  $P^{\Gamma}$ are all different. It follows
that $C_t$ defines a smoothing of $(C, 0)$ by Theorem
\ref{vlocal}.

In  \cite{VirL} Viro introduces sqh-smoothing as follows: Suppose
that $Q = \sum_{(i,j) \in \D^- \cap \Z^2} c_{i,j} x^i y^j +
P_0^{\Gamma} \in \R[x, y]$, is non degenerated with respect to its
Newton polygon $\D$. Let  $\w: \D \cap \Z^2 \rightarrow \Z$ be the
function:
\begin{equation} \label{omega}
  \w
 (r, s) = nm - n r - m s.
\end{equation}
Then $P_t := \sum_{(i,j) \in \D^- \cap \Z^2} c_{i,j} t^{\w(i,j)}
x^i y^j + P_0(x, y)$ defines a sqh-smoothing of $(C,0)$. For
technical reasons we consider in Definition \ref{sqh} a slightly
more general notion, by allowing terms which depend on $t$ and
which have exponents above the lower part of the Newton polyhedron
of $P_t$, viewed in $\R^3$
 (cf. Notations \ref{notation-viro}).

\subsection{Size of ovals of sqh-smoothings}

In this Section we consider some remarks about the sizes of ovals
of sqh-smoothing $C_t$ of $(C,0)$ We say that $a(t) \sim t^{\g}$
if there exists a non zero constant $c$ such that $a(t) \sim c
t^{\g}$ when $t>0$ tends to $0$.
\begin{Defi} \label{size-defi}
An oval of $C_t$ is of size $(t^\a, t^\b)$ if it is contained in a
minimal box of edges parallel to the coordinate axis such that
each vertex of the box has coordinates of the form $(\sim t^{\a},
\sim t^{\b})$.
\end{Defi}

\begin{Pro} \label{size-sqh}
Let  $C_t$ is a sqh-smoothing of $(C,0)$.  If the sqh-smoothing is
described with Notations \ref{notation-sqh}, then each oval of the
smoothing is of size $(t^n, t^m)$.
\end{Pro}
{\em Proof.} The critical points of the projection $C_t
\rightarrow \C$, given by $(x, y) \mapsto x$, are those defined by
$P_t =0$ and $\frac{\partial P_t}{\partial y} =0$. The critical
values of this projection are defined by zeroes in $x$ of the
discriminant $\D_y P_t$. Using the non degeneracy conditions on
the edges of the local Newton polyhedron of $P_t (x, y)$, viewed
as a polynomial in $x,y, t$ and Théorème 4 in \cite{GP00} we
deduce that the local Newton polygon of $\D_y P_t$, as a
polynomial in $x, t$ has only two vertices: $((n-1)m, 0)$ and $(0,
(n-1) nm)$. It follows by the Newton-Puiseux Theorem that the
roots of the discriminant $\D_y P_t$ as a polynomial in $x$,
express as fractional power series of the form:
\begin{equation} \label{aa}
x =  t^n  \epsilon_r, \mbox{ where } \epsilon_r\in \C \{ t^{1/k}
\} \mbox{ and } \epsilon_r (0) \ne 0 , \mbox{ for } r =1, \dots ,
(n-1) m,
\end{equation}
which correspond to the root $x=0$ of $\D_y P_0$ and
\begin{equation} \label{bb}
x = \varepsilon_s  \in \C \{ t^{1/k} \} \mbox{ with }
\varepsilon_s (0) \ne 0
\end{equation}
which correspond to the non zero roots $\varepsilon_s (0)$  of
$\D_y P_0$ (counted with multiplicity). The critical values of the
smoothing are among those described by (\ref{aa}). The critical
values (\ref{bb}) correspond to slight perturbations of critical
values of $C = C_0 \rightarrow \C$, outside the Milnor ball of
$(C,0)$.

We argue in a similar manner for the projection $C_t \rightarrow
\C$, given by $(x,y) \mapsto y$. $\Box$

\section{Harnack smoothings} \label{harnack-smoothings}

We generalize the notions of {\em maximal position} and {\em good
oscillation} and {\em Harnack}, introduced in Section \ref{Mik},
to the case of smoothings as follows:

\begin{Defi} \label{hsm}
\begin{enumerate}
\item[(i)] A smoothing $C_t$ of $(C,0)$ is in  {\em maximal
position (resp. has good oscillation) with respect to a line} $L$
passing through the origin if there exists $(L,C)_0$ different
points of intersection of $C_t $ with $L$, which tend to $0$ as
the parameter $t>0 $ tends to $0$,  and which are all contained in
an arc $\mathbf{a} \subset \R C_t$ of the smoothing (resp. $C_t$
is in maximal position and the order of the points in $\R L \cap
\R C_t$ is the same when the points are viewed in the line $L$ or
in the arc $\mathbf{a}$), for all $0 < t  \ll 1$.

\item[(ii)] A smoothing $C_t$ of $(C,0)$ is in  {\em maximal
position with respect to two lines} $L_1, L_2$ passing through the
origin if it is in maximal position with respect to the lines
$L_1$ and $L_2$ and if the $(L_i,C)_0$  points of intersection of
$C_t $ with $L_i$, which tend to $0$ as the parameter $0 <t $
tends to $0$, are all contained in an arc $\mathbf{a}_{i,t}$ of
the smoothing $\R C_t$ , for $i=1, 2$, such that
$\mathbf{a}_{1,t}$ and $\mathbf{a}_{2,t}$ are disjoint and
contained in the same component of the smoothing $\R C_t$, for all
$0 < t \ll 1$.

\item[(iii)] A Harnack smoothing is a $M$-smoothing which is in
maximal position with good oscillation with respect to the lines
$L_1$ and $L_2$.
\end{enumerate}
\end{Defi}

\begin{figure}[htbp]
$$\epsfig{file=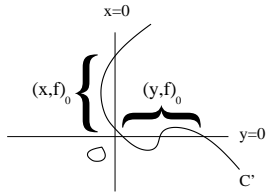, height= 25 mm}$$
\caption{A Harnack smoothing of the cusp $ y^2 - x^3 =0 $
\label{hk}}
\end{figure}

\begin{Rem}
 Every real
plane branch admits a Harnack smoothing (this result is implicit
in the blow up method of \cite{RisH}).
\end{Rem}

\subsection{The case of non degenerated sqh-polynomials with peripheral roots of the same sign}

We consider the case of a plane curve singularity $(C,0)$ defined
by a semi-quasi-homogeneous polynomial $P_0 \in \R[x, y]$ such
that the peripheral roots  associated to the compact edge of its
local Newton polygon are all real, different and of the same sign.
We prove that there exists a Harnack smoothing $C_t$ of $(C,0)$
constructed by Patchworking. This result is quite similar to
\cite{KRS} Theorem 4.1 (1), where under the same hypothesis, they
prove that a $M$-smoothing exists. We prove  then that the
embedded topological type of the Harnack smoothing of $(C,0)$ is
unique.

We keep Notations \ref{notation-sqh} in the description of the
polynomial $P_0$.

\subsubsection{Existence of Harnack smoothings}

\begin{figure}[htbp]
$$\epsfig{file=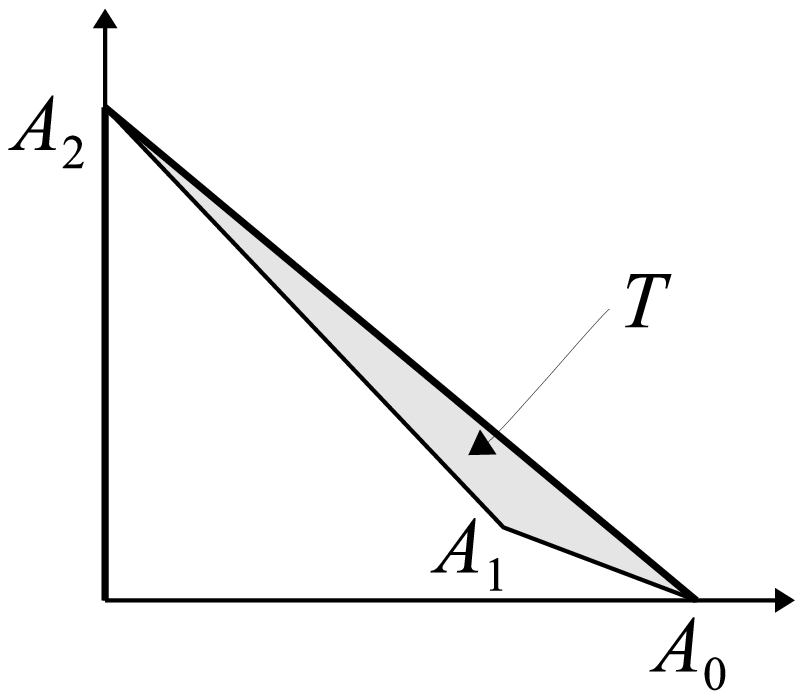, height= 25 mm}$$
\caption{\label{triangle}}
 \end{figure}

\begin{Lem} \label{convex}
There exists a piece-wise affine linear convex function $\w: \D
\rightarrow \Z_{\geq 0}$, which takes integral values on $\D \cap
\Z^2$, vanishes on $\Gamma \cap \Z^2$ and induces a triangulation
of $\D$ with the following properties:
\begin{enumerate}
\item[(i)] All the integral points in $\D^-$ are vertices of the
triangulation.

\item[(ii)] There exists exactly one triangle $T$ in the
triangulation which contains $\Gamma$ as an edge. The triangle $T$
is transformed by a translation and a $SL(2, \Z)$-transformation
into the triangle $[(0,1),(e,0),(0,0)]$.

\item[(iii)] If $T' \ne T$ is in the triangulation then $T'$ is
primitive.
\end{enumerate}
\end{Lem}
{\em Proof.} With the above notations take $A_1 \in \D^-$ the
closest  integral point to $\Gamma = A_0 A_2$. Let $T$ be the
triangle with vertices $A_0,A_1, A_2$ (see Figure \ref{triangle}).
Then assertion (ii) follows. It is then easy to construct a convex
triangulation of $\D$ which contains $T$ and which is primitive on
$\D \setminus T$ (see \cite{KRS}).   $\Box$

We say that a distribution of signs $\epsilon$ is {\em compatible}
with a polynomial $Q = \sum c_{i,j} x^i y^j \in \R [x, y]$ if
$\mbox{{\rm sign}} (c_{i,j}) = \epsilon(i,j)$.

\begin{Pro} \label{compatible} With the above notations,
let  $P_0 \in \R[x, y]$ be a semi-quasi-homogeneous polynomial
defining a plane curve singularity $(C,0)$. If the peripheral
roots of $P^{ \Gamma}$ are all different and of the same sign,
then there are two Harnack distribution of signs, $\epsilon_{1}$
and $\epsilon_2$, which are compatible with $P^{ \Gamma}$. Let
$\w$ be as in Lemma \ref{convex}. Then the polynomial
\[
P_t = \sum_{(i, j) \in \D^- \cap \Z^2}  \epsilon_k (i,j)  t^{\w(i,
j)} x^i y^j + P_{0},
\]
defines a Harnack smoothing of $(C,0)$, for $k=1,2$.
\end{Pro}
{\em Proof.} We keep Notations \ref{notation-sqh}. We suppose
without loss of generality that the peripheral roots of $P^\Gamma$
are all positive.  Then the signs of the coefficients
corresponding to consecutive terms in the edge $\Gamma$ are always
different. By a simple observation on the set of Harnack
distribution of signs  it follows that there exists precisely two
different  Harnack distribution of signs $\epsilon, \epsilon'$
which are compatible with $P^\Gamma$ (see Definition
\ref{harnack-sign}). Since $P^{\Gamma} (x, y) = P^\Gamma \circ
\r_\Gamma (x, y)$ by definition of $\r_\Gamma$ (see Section
\ref{toric}) the distributions of signs $\epsilon$ and $\epsilon'$
are related by $\epsilon' = \epsilon \circ \r_\Gamma$.

We use notations \ref{notation-viro}. Let  $\w : \D \rightarrow
\R_{\geq 0}$ be a piece-wise affine  convex  function satisfying
the statement of Lemma \ref{convex}. By Lemma \ref{convex} the
chart of $P^{\hat{T}}_{{1}}$ is transformed by a translation and a
$SL(2, \Z)$-transformation to the chart of a polynomial with
Newton polygon with vertices $(0,0)$, $(0,1)$ and $(e,0)$, i.e.,
to the chart of the  graph of a polynomial of one variable with
$e$ different positive real roots . The topology of this chart is
determined by the sign of the term corresponding to $(0,1)$. Let
${\w}':\D \rightarrow \R$ be a convex piece-wise affine function
determined by its integral values on $\D \cap \Z^2$. We can assume
that $\w'$ induces a primitive triangulation of $T$ and of $\D
\setminus \mbox{\rm int} T$ the functions $\w$ and ${\w}'$ define
the same primitive triangulation (by translating $\w$ we can
assume that $\w_{| T} =\a > 0$, then define $\w'$ as $\w'= \w$ on
$\D \setminus T$ and $\w'_{| A_0 A_2}$ be strictly convex piece
wise-linear and positive).

 The topology of the chart of
$P^{\hat{T}}_{{1}}$  coincides with the topology of the chart of
the polynomial $\sum_{(i, j) \in T \cap \Z^2} \epsilon_k (i,j)
t^{\tilde{\w}(i, j)} x^i y^j$.
 Notice that the patchworking
of $P_t$ is combinatorial for all triangles $T'$ of the
subdivision with the possible exception of $T$. It follows from
Theorem \ref{vlocal} that the signed topological type of the
smoothing $P_t$ coincides with that of the smoothing defined by
\[
Q_t := \sum_{(i, j) \in \D \cap \Z^2} \epsilon_k (i,j)
t^{\tilde{\w}(i, j)} x^i y^j,
\]
which is constructed by combinatorial patchworking. It follows
that the polynomial $P^{\hat{\D}}_{{1}}$ defines a Harnack curve
in the toric surface $Z(\D)$,  by Proposition \ref{harnack-curve}.
Therefore $P_t$ defines then a Harnack smoothing of $(C,0)$ with
$\# \mbox{{\rm int}} \D \cap \Z^2$ ovals and $e$ non compact
components. $\Box$

\subsubsection{The topological type of a Harnack
smoothing}

\begin{The} \label{1par}
Let $(C,0)$ be a plane curve singularity defined by a
semi-quasi-homogeneous polynomial $P_0 \in \R[x, y]$ non
degenerated with respect to its local Newton polygon. We denote by
$\Gamma$ the compact edge of this polygon. We suppose that the
peripheral roots of $P_0^{\Gamma}$ are all real. Let $C_t$ define
a semi-quasi-homogeneous $M$-smoothing of $(C,0)$ such that $C_t$
is in maximal position with respect to the coordinate lines. If
$B$ denotes a Milnor ball for $(C,0)$ then we have that:
\begin{enumerate}
\item[(i)] The peripheral roots of $P_0^{\Gamma}$ are of the same
sign.

 \item[(ii)] The polynomial $P^{\hat{\D}}_{{1}}$ defines a Harnack
curve in $Z(\D)$.

\item[(iii)] The smoothing $C_t$ is Harnack.

\item[(iv)] There is a unique topological type of triples $(\R
{B}, \R C_t , B \cap (\R^*)^2)$.

\item[(v)] The topological type of the smoothing $C_t$ is
determined by $\D$.

\end{enumerate}
\end{The}
{\em Proof.} We follow Notations \ref{notation-sqh} to describe
the polynomial $P_0$.
Since the peripheral roots of the polynomial $P_0^{\Gamma}$ are
all real it follows that the singularity $(C,0)$ has exactly $e$
analytic branches which are all real. It follows that the
smoothing $C_t$ has $e$ non closed components. If $C_t$ is
$M$-smoothing there are precisely $\# (\mbox{\rm int}{\D} \cap
\Z^2)$ ovals by (\ref{mu-delta}). Since the smoothing $C_t$ is
Harnack, i.e., it is in maximal position with respect to the
coordinate axis, none of these ovals cuts the coordinate axis.

We consider the curve $\tilde{C}$ defined by  the polynomial
$P^{\hat{\D}}_{{1}}$ in the real toric surface $Z (\D)$. By
Theorem \ref{vlocal} and Remark \ref{viro-toric} there are exactly
$\# (\mbox{\rm int}{\D} \cap \Z^2)$ ovals in the chart
$\mbox{Ch}_\D ^* (\tilde{C})$. These ovals, when viewed in the
toric surface $Z(\D)$ by Proposition \ref{quatro}, do not meet any
of the coordinate lines of $Z(\D)$. By the same argument we have
that the curve $\tilde{C}$
 is in maximal position with respect to
the two coordinate lines $x=0$ and $y=0$ of $Z(\D)$, corresponding
respectively to the vertical and horizontal edges of $\D$. It
follows that the number of components of $\R \tilde{C}$ is $\geq 1
+ \# (\mbox{\rm int}{\D} \cap \Z^2)$, which is equal to the
maximal number of components, hence $\tilde{C}$ is a $M$-curve in
$Z(\D)$.

We deduce that the non compact connected components of $\R
\tilde{C} \cap (\R^*)^2 \subset Z(\D)$ glue up in one connected
component of $\R \tilde{C}$. This component contains all the
intersection points with the coordinate lines of $Z(\D)$, since
$\tilde{C}$ is in maximal position with respect to the lines $x=0$
and $y=0$, and by assumption the peripheral roots of $P_0^\Gamma$
are all real.

By definition of maximal position with respect to two lines, there
is exactly one component of the smoothing $C_t$  containing all
the points of intersection of $C_t$ with the coordinate axis. We
deduce the following assertion by translating this information in
terms of the chart of $\tilde{C}$, using  Theorem \ref{vlocal} and
Proposition \ref{quatro}:  there are two disjoint arcs
$\mathbf{a}_x$ and $\mathbf{a}_y$ in $\R \tilde{C}$ containing
respectively the points of intersection of $\tilde{C}$ with the
toric axis $x=0$ and $y=0$, which do not contain any point in $\R
Z (\Gamma)$ (otherwise there would be more than one non compact
component of the smoothing $C_t$ intersecting the coordinate axis,
contrary to the assumption of maximal position). It follows from
this that the curve $\tilde{C}$ is in maximal position with
respect to the coordinate axis in $Z (\D)$. Mikhalkin's Theorem
\ref{Mikhalkin} implies the second assertion. Then the other three
assertions are deduced from this by Theorem \ref{vlocal} and
Proposition \ref{top-harnack}. $\Box$

\begin{figure}[htbp]
 $$\epsfig{file=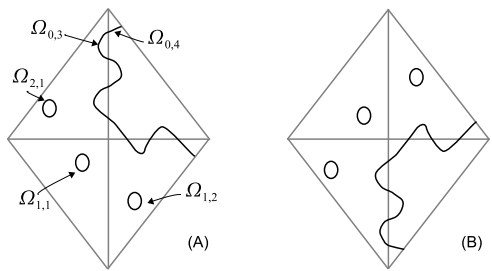, height= 35 mm}$$
 \caption{Charts of Harnack smoothings of $y^4 - x^3=0$.
 The signed topological  type in (A) is normalized      \label{normalizado}}
 \end{figure}

\begin{Rem} \label{top-harnack2} With the hypothesis and notations of Theorem \ref{1par},
we have that the connected components $\O_{i,j}$, for $(i, j) \in
\D \cap \Z^2$, of chart $\mbox{Ch}_\D ^* (\tilde{C})$,
are described by Proposition \ref{top-harnack}. If the peripheral
roots of $P_0^\Gamma$ are positive,  up to replacing $P^{\D}_t$ by
$P^{\D}_t \circ \r_{\Gamma}$, one can always have that $\O_{0,n}
\subset \R^2_{0,0}$ and then:
\begin{equation} \label{indice-normalized-smoothing}
\O_{r,s} \subset   \R^2_{n+s,r}, \quad \forall ( r,s) \in \D \cap
\Z^2.
\end{equation}
Otherwise, we have that:
\begin{equation} \label{indice-smoothing}
\O_{r,s} \subset \r_\Gamma (  \R^2_{n+s,r}), \quad \forall ( r,s)
\in \D \cap \Z^2.
\end{equation}
Compare in Figure \ref{normalizado} the position of $\O_{0,4}$ in
(A) and (B).
\end{Rem}

\begin{Defi} \label{normalized-smoothing}
If (\ref{indice-normalized-smoothing}) holds then we say that the
signed topological type of the smoothing $C_t$ (or of the chart of
$P_{\D}$),  is {\em normalized} (cf. Definition
\ref{topological2}).
\end{Defi}

As an immediate corollary of Theorem \ref{1par} we deduce that:
\begin{Pro}  \label{signed-topological}
There are two signed topological types of sqh-$M$-smoothings of
$(C,0)$ in maximal position with respect to the coordinate axis.
These types are related by the orthogonal symmetry $\r_\Gamma$ and
only one of them is normalized.
\end{Pro}

One of the aims of this paper is to study to which extent Theorem
\ref{1par} admits a valid formulation in the class of real plane
branch singularities. In general the singularities of this class
are degenerated with respect to their local Newton polygon, in
particular we cannot apply Viro's method to those cases.
Classically smoothing of this type of singularities is constructed
using the blow up construction. We present in the following
sections an alternative method which applies Viro method at a
sequence of certain infinitely near points.

\section{A reminder on toric resolutions of real plane
branches}\label{toric-res}

We recall the construction of an embedded resolution of
singularities of a real plane branch  by a sequence of local toric
modifications. For details see \cite{AC} and  \cite{roots}. See
\cite{Oda}, \cite{Fulton} for more on toric geometry and
\cite{Oka1}, \cite{Oka2}, \cite{LO} and \cite{Rebeca} and for more
on toric geometry and plane curve singularities.

A germ $(C, 0)$ of real plane curve, defined by $F =0$ for $F \in
\R[x, y]$, defines a real plane branch if it is analytically
irreducible in $(\C^2, 0)$ and if it admits a real Newton Puiseux
parametrization (normalization map):
\begin{equation} \label{parametrization}
\left\{
\begin{array}{lcl}
x (t) & = & t^{e_0},
\\
y  (t) & = & \sum_{i} \eta_i t^i, \quad \mbox{ with } \eta_i \in
\R.
\end{array}
\right.
\end{equation}

If the coordinate line $x=0$ is not tangent to $(C,0)$ then $e_0 =
\mbox{ord}_t (x(t))$ is the {\em multiplicity} of $(C,0)$. By
(\ref{parametrization}) and a suitable change of coordinates, we
have that $C$ has an equation $F =0$, with
\begin{equation} \label{normal1}
F = ( y ^{n_1} - x^{m_1})^{e_1} + \cdots,
\end{equation}
such that $\mbox{\rm gcd } (n_1, m_1)=1$, the integer $e_0 := e_1
n_1$ is the intersection multiplicity with the line, $x=0$, and
the terms which are not written have exponents $(i, j)$ such that
$i n_1 + j m_1
> n_1 m_1 e_1$, i.e., they lie above the compact edge
\[
\Gamma_1 := [(0, n_1 e_1), (m_1 e_1, 0)]
\] of the local Newton
polygon  of $F$.

\begin{figure}[htbp]
$$\epsfig{file=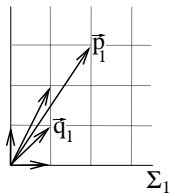, height= 25 mm}$$
\caption{The subdivision $\Sigma_1$ associated to  $F$ in Example
\ref{10}  \label{sigma2-3}}
\end{figure}

The vector $\vec{p}_1 = (n_1, m_1)$ is  orthogonal to $\Gamma_1$
and defines a subdivision of the  positive quadrant $\R^2_{\geq
0}$, which is obtained by adding the ray $\vec{ p}_1 \R_{\geq 0}$.
The quadrant $\R^2_{\geq 0}$ is subdivided in two cones, $\t_i
:=\vec{ e}_i \R_{\geq 0} + \vec{ p}_1 \R_{\geq 0}$, for $i=1, 2$
and $\vec{ e}_1, \vec{ e}_2$ the canonical basis of $\Z^2$.  We
define the {\em minimal regular subdivision} $\Sigma_1$ of
$\R^2_{\geq 0}$ which contains the ray $\vec{ p}_1 \R_{\geq 0}$ by
adding the rays defined by those integral vectors in $\R^2_{>0 }$,
which belong to the boundary of the convex hull of the sets $(\t_i
\cap \Z^2) \backslash \{ 0\}$, for $i=1,2$. We denote by $\s_1$
 the unique
cone of $\Sigma_1$ of the form, $\s_1 : = \vec{p}_1 \R_{\geq 0} +
\vec{q}_1 \R_{\geq 0}$ where $\vec{q}_1 = (c_1, d_1)$ satisfies
that:
\begin{equation} \label{ave-i}
c_1 m_1 - d_1 n_1 = 1.
\end{equation}
See an example in Figure \ref{sigma2-3}.

By convenience we denote $\C^2$ by $Z_1$.
 We define in the sequel
a sequence of proper birational maps $\Pi_j: Z_{j+1} \rightarrow
Z_j$, for $j=1, \dots, g$ such that the composition $\Pi_g \circ
\dots \circ \Pi_1$ is an {\em embedded resolution} of the germ
$(C,0)$, i.e., the pull back of the germ $(C,0)$ by  $\Pi_g \circ
\dots \Pi_1$ is a normal crossing divisor in the smooth surface
$Z_{g+1}$.

In order to describe these maps we denote the coordinates $(x, y)$
by $(x_1, y_1)$ and the origin $0 \in \C^2 = \Z_1$ by  $o_1$. We
also denote $F$ by $F^{(1)}$ and $C$ by $C^{(1)}$.  The
subdivision $\Sigma_1$ defines a proper birational map $\Pi_1: Z_2
\rightarrow Z_1$, which is obtained by gluing maps $\C^2 (\s)
\rightarrow \C^2$, where $\s$ runs through the set of two
dimensional cones in $\Sigma_1$. For instance, the map $\p_1: \C^2
(\s_1) \rightarrow \C^2$ is defined by
\begin{equation} \label{mm-i}
\begin{array}{lcl}
x_1 & = & u_2 ^{c_1} x_2^{n_1},
\\
y_1 & =  & u_2^{d_1} x_2^{m_1},
\end{array}
\end{equation}
where $u_2, x_2$ are canonical coordinates for the affine space
$\C^2 (\s_1)$.

It should be noticed that the map $\Pi_1$ is a composition of
point {\em blow-ups}, as many as rays added in $\Sigma_1$ to
subdivide $\R^2_{\geq 0}$. Each ray  $ \vec{a}\R_{\geq 0} \in
\Sigma_1$ corresponds bijectively to a projective line $\C P^1$,
embedded as an irreducible component of the exceptional divisor
$\Pi_1^{-1} (0)$. We denote by $E_2 \subset Z_2$ the exceptional
divisor defined by $x_2=0$  in the chart $\C^2 (\s_1)$, the other
exceptional divisor in this chart being defined by $u_2 =0$.
Notice that the point at the infinity of the line $x_2 =0$, for
instance, is the origin of the chart $\C^2 (\s_1 ')$, where $\s_1'
\in \Sigma_1$ is the two dimensional cone adjacent to $\s_1$ along
the ray $\vec{p}_1 \R_{\geq 0}$.

We have that $\Pi_1 ^* (C^{(1)})$ defines a {\em Cartier divisor}
on $Z_2$. For instance, on $\C^2 (\s_1)$  it is defined by $F
\circ \p_1 =0$. The term  $F \circ \p_1 =0$ decomposes as:
\begin{equation} \label{strict-j}
F^{(1)} \circ \p_1 =   \mbox{\rm Exc} (F^{(1)}, \p_1)
\bar{F}^{(2)} (x_2, u_2), \mbox{ where } \bar{F}^{(2)} (0,0) \ne
0,
\end{equation}
and
\[ \mbox{\rm Exc} (F^{(1)}, \p_1) := y_1^{e_0} \circ \p_1.
\]
The polynomial $\bar{F}^{(2)} (x_2, u_2)$ defines the {\em strict
transform} $C^{(2)}$ of $C^{(1)}$, i.e.,  the closure of the
pre-image by $\p_1^{-1}$ of the punctured curve $C^{(1)} \setminus
\{ 0 \}$ on the chart $\C^2 (\s_1)$.  The function $ \mbox{\rm
Exc} (F^{(1)}, \p_1) $ defines the exceptional divisor of $\Pi_1^*
(C^{(1)})$ on this chart.

We analyze the intersection of the strict transform with the
exceptional divisor on the chart $\C^2 (\s_1)$: using
(\ref{normal1}) we find that $ \bar{F}^{(2)} (x_2, 0) =1$ and
\[
 \bar{F}^{(2)}  ( 0, u_2) = ( 1  - u_2^{c_1 m_1 - d_1 n_1}) ^{e_1}
\stackrel{\mbox{ (\ref{ave-i})}}{=}  (1-u_2) ^{e_1}.
\]
By a similar argument on the other charts it follows that $E_2$ is
the only exceptional divisor of $\Pi_1$ which intersects the
strict transform $C^{(2)}$ of $C$, precisely at the point $o_2$ of
the chart $\C^2 (\s_1)$ with coordinates $x_2 =0$ and $u_2=1$,
with intersection multiplicity equal to $e_1$. If $e_1 =1$ then
the strict transform is smooth at $o_2$ and the intersection with
the exceptional divisor $x_2 =0$ is transversal, hence the divisor
$\Pi_1 ^* (C)$ has smooth components which intersect
transversally. In this case, the map $\Pi_1$ is an {\em embedded
resolution} of the germ $(C,0)$ by definition.

We define a pair of real coordinates $(x_2, y_2)$ at the point
$o_2$, where
\begin{equation} \label{xi-rel}
\xi_2 y_2 = 1- u_2 + x_2 R_2 (x_2, u_2)  \mbox{ for some } \xi_2
\in \R^* \mbox{ and } R_2 \in \R [x_2, u_2],
\end{equation}
such that $C^{(2)}$ is defined by a polynomial, which we call {\em
the strict transform function}, of the form:
\begin{equation}\label{fg}
 F^{(2)} (x_2, y_2) = (y_2^{n_2} -
x_2^{m_2} ) ^{e_2} +  \cdots,
\end{equation}
where $\mbox{\rm gcd} (n_2, m_2)=1$,  $e_1 = e_2 n_2$  and the
terms which are not written have exponents $(i, j)$ such that $i
n_2 + j m_2
> n_2 m_2 e_1$, i.e., they lie above the compact edge
$\Gamma_2 $ of the (local) Newton polygon of $F^{(2)}(x_2, y_2)$.
The result of substituting in  $F^{(2)} (x_2, y_2)$, the term
$y_2$ by  using (\ref{xi-rel}), is equal to $\bar{F}^{(2)}(x_2,
u_2)$.

We can iterate this procedure defining for $j>2$ a sequence of
toric maps $\Pi_{j-1} :Z_{j} \rightarrow Z_{j-1}$, which are
described by replacing the index $1$ by $j-1$ and the index $2$ by
$j$ above. In particular, when we refer to a Formula, like
(\ref{ave-i}) at level $j$, we mean after making this replacement.

We denote by $\mbox{\rm Exc} (F^{(1)} , \p_1 \circ \cdots \circ
\p_{j}
 ) $ the {\em exceptional function} defining the exceptional
 divisor of $(\Pi_1 \circ \cdots \circ \Pi_{j})^* (C)$ on the
 chart $\C^2 (\s_{j}) \subset Z_j$. We have that
 \begin{equation} \label{EJ}
 \mbox{\rm Exc} (F^{(1)}, \p_1 \circ \cdots \circ \p_{j} ) =  (y_1^{e_0}
 \circ \p_{1} \circ \cdots \circ \p_{j}) \cdots (y_{j}^{e_{j-1}}
 \circ \p_{j} ).
 \end{equation}

Since by construction we have that $e_j | e_{j-1} | \cdots | e_1
|e_0$ (for $ | $ denoting divides), at some step we reach a first
integer $g$ such that $e_g =1$ and then the process stops. The
composition of blow ups $\Pi_g \circ \dots \circ \Pi_1$ is an {\em
embedded resolution} of the germ $(C,0)$. It is {\em minimal}, in
the number of exceptional divisors required, if one of the
coordinate axis $x_1 =0$ or $y_1 =0$ is not tangent to $(C,0)$, in
particular if $e_0$ is the multiplicity of $(C,0)$ (see
\cite{AC}).

If $x_1=0$ is not tangent to $C$ then the sequence of pairs
$\{(n_j, m_j)\}_{j=1}^g$ classifies the embedded topological type
of the pair $(C,0) \subset (\C^2, 0)$. This is because the
$\{(n_j, m_j)\}_{j=1}^g$  determine and are determined by the
classical characteristic pairs of a plane branch (see \cite{AC}
and  \cite{Oka1}).

\begin{Exam} \label{10}
The embedded resolution of the real plane branch singularity $(C,
0)$ defined by $F= (y_1^2 -x_1^3 )^3 - x_1^{10} =0$ is as follows.
\end{Exam}
The morphism $\p_1$ of the toric resolution is defined by
\[
\begin{array}{lcl}
x_1 & = & u_2 ^{1} x_2^{2},
\\
y_1   & =  & u_2^{1} x_2^{3}.
\end{array}
\]
Let $z_1 := y_1^2 - x_1^3 $, then we have $z_1 \circ \p_1 = u_2 ^2
x_2^ 6 ( 1 - u_2) =u_2 ^2 x_2 ^ 6 y_2$, where $y_2 := 1 -u_2$
defines the strict transform function of $z_1$, and together with
$x_2$ defines local coordinates at the point of intersection $o_2$
with the exceptional divisor $x_2=0$.

For $F$ we find that:
\[
\begin{array}{lcl}
F \circ \p_1 &  = & u_2 ^6 x_2 ^{18 } \left( (1 - u_2)^3 - u_2^{4}
x_2 ^{2} \right).
\end{array}
\]
Hence  $\mbox{\rm Exc} (F, \p_1) :=  y_1^6 \circ \p_1 = u_2 ^6 x_2
^{18 }$ is the exceptional function associated to $F$, and
\[ F ^{(2) } = y_2^3 - (1- y_2)^4 x_2^2  \] is the strict
transform function.

\begin{Not} \label{notation-res}
For  $j=1, \dots,g$:
\begin{enumerate}
\item[(i)] Let $ \Gamma_j = [(m_{j} e_{j}, 0) , ( 0, n_j e_j)]$,
be the unique compact edge of the local Newton polygon of $F^{(j)}
(x_j, y_j)$ (see (\ref{fg}) at level $j$).

\item[(ii)] Let $\D_j$ the Newton diagram of $F^{(j)} (x_j, y_j)$.
We denote by  $\D_j^- $ the set $\D_j^- = \D_j \setminus
\Gamma_j$.

\item[(iii)] Let $\Xi_j = [(0,0), (0, n_j e_j)]$ be the edge of
$\D_j$ which is the intersection of $\D_j$ with the vertical axis.

\item[(iv)]
 Let $\w_j: \D_{j} \cap \Z^2 \rightarrow \Z$ be defined by
 \begin{equation} \label{omega-j}
 \w_j
(r, s) = e_{j}  ( e_{j} n_{j} m_{j} - r n_{j} - s m_{j} )
\end{equation}
\end{enumerate}
\end{Not}
Notice that (\ref{omega-j}) is defined analogously as
(\ref{omega}).

 \begin{figure}[htbp]
$$\epsfig{file=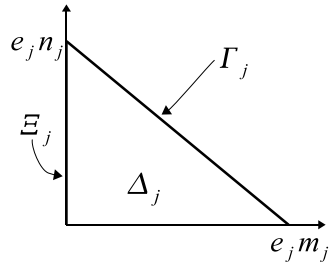, height= 35 mm}$$
\caption{The Newton diagrams in Notation \ref{notation-res}}
 \end{figure}

\begin{Pro}
The following formula for the Milnor number of $(C,0)$ at the
origin is deduced in \cite{roots}.
\begin{equation} \label{mu-delta-j}
\mu (C)_{0} =  2 \sum_{j=1}^{g}   \# ( \mbox{{\rm int}} \D_j \cap
\Z^2 ) + e_{j} -1.
\end{equation}
\end{Pro}

\subsection{A set of polynomials defined from the embedded resolution}

We associate in this section some polynomials to the elements in
$\D_j^{-} \cap \Z^2$. From these polynomials we define a class of
deformations which we will study in the following sections.

\begin{Lem} \label{comb2}
If $(r, s) \in \Z^2_{\geq 0}$ with $s < e_{j-1}$
there exist $M_j (r, s) \in \R[x_1, y_1]$ and integers $k_2>0,
\dots , k_j  >0$ such that:
\begin{equation} \label{ler}
\mbox{\rm Exc} (F, \p_1 \circ \cdots \circ \p_{j-1} )  \,
u_2^{k_2} \cdots u_j^{k_j} \, x_j^r \, y_j^s = M_j (r,s) \circ
\pi_1 \circ \cdots \circ \pi_{j-1}.
\end{equation}
\end{Lem}

To avoid cumbersome notations  we denote simply by $u_i$ the term
$u_i \circ \p_{i} \circ \cdots \circ \p_{j-1}$, whenever $j > i $
is clear form the context. For instance we have done this in
Formula (\ref{ler}) above. In particular, by (\ref{xi-rel}) at
level $\leq j$, the restriction of the function $u_i$ to $x_j =0$
is equal to:
\begin{equation} \label{units}
u_i = \left\{
\begin{array}{lcl}
1 & \mbox{ if } &  i < j,
\\
1 - \xi_j y_j &  \mbox{ if } &  i=j
\end{array}
\right.
\end{equation}

\begin{Rem} \label{dep}
The integers $k_2, \dots, k_j$ depend on $(r,s)$ and on the
singular type of the branch $(C,0)$. They can be determined
algorithmically, and are unique up to certain conditions on the
polynomials $M_j (r,s)$. The polynomials $M_j(r,s)$ are
constructed as monomials in $x_1, y_1$ and some Weierstrass
polynomials defining {\em curvettes} at certain irreducible
exceptional divisors of the embedded resolution of $(C,0)$ (see
\cite{roots}). For more details on the construction of these
curves and their applications see \cite{PP}, \cite{Rebeca},
\cite{Z} and \cite{Moh}.
\end{Rem}

\begin{Exam} \label{11}
The following table indicates the terms $M_2(r,s)$ for $(r,s) \in
\D_2$ corresponding to Example \ref{10}. The symbol $z_1$ denotes
the first approximate root $y_1^2 - x_1^3$.
\end{Exam}

\begin{center}
\begin{tabular}{|c|c|c|c|c|c|}
\hline
  $(r,s)$ &   $(0,0) $ & $(0,1)$ & $(0,2) $& $(1,1)$ &$ (1,0)$
\\
\hline $ M_2 (r,s)$  & $x_1^9$ & $x_1^6 z_1 $& $x_1^3 z_1^2 $&
$x_1^5 y_1 z_1$ & $x_1^8 y_1$
\\
\hline
\end{tabular}
\end{center}
For instance, we have that $M_2 (1, 1)=  x_1^5 y_1 z_1$, since
$x_1^5 y_1 z_1 \circ \p_1 = \mbox{\rm Exc} (F^{(1)}, \p_1)  u_2^2
x_2 y_2$, where $\mbox{\rm Exc} (F^{(1)} , \p_1) =  u_2 ^6 x_2
^{18 }$ by Example \ref{10}.

\section{Multi-semi-quasi-homogeneous smoothings of a real plane
branch} \label{last}

In this Section we introduce a class of deformations of a plane
branch $(C,0)$, called multi-semi-quasi-homogeneous (msqh)
deformations and we describe their basic properties.

We suppose that $(C,0)$ is a plane branch  defined by an equation,
$F(x, y) =0$, such that its embedded resolution consists of $g$
toric maps (see Section \ref{toric-res} and Notation
\ref{notation-res}). Consider the following algebraic expressions
in terms of the polynomials of Lemma \ref{comb2} as a sequence of
deformations of the polynomial $F(x_1, y_1)$ defining $C$. We
denote by $M_1(r,s)$ the monomial $x_1^r y_1^s$ and $\tr_j$ will
denote $t_j, \dots, t_g$ for any $1 \leq j \leq g$.
\begin{equation} \label{pg}
\left\{
\begin{array}{ccccccl}
P_{\tr_g} & := & F& +  & \displaystyle{\sum_{(r,s) \in \D_{g}^-
\cap \Z^2}} \, A_{r,s}^{(g)}  \, & t_{g}^{\w_{g} (r,s)} \,
 & M_{g} (r,s)
\\
 P_{\tr_{g-1}} & := & P_{\tr_g} & +  & \displaystyle{\sum_{(r,s) \in \D_{g-1}^-
\cap \Z^2}} \, A_{r,s}^{(g-1)}  \, & t_{g-1}^{\w_{g-1} (r,s)} \,
_{0\leq s < e_{g-1}}^{0 \leq r} \a_{g-1, r, s} (\l_{g-1})  \,
 & M_{g-1} (r,s)
\\
\dots & \dots & \dots & \dots & \dots & \dots & \dots
\\
 P_{\tr_1 } & := &  P_{\tr_2} & +  & \displaystyle{\sum_{(r,s) \in \D_{1}^- \cap
\Z^2}} \,  A_{r,s}^{(1)} \, &  t_{1}^{\w_{1} (r,s)} \,
& M_1(r,s).
\end{array}
\right.
\end{equation}
The terms $A_{r,s}^{(j)}$ are some real constants while the $t_j$
are real parameters, for $j=1, \dots, g$. For technical reasons we
will suppose that $A_{0,0}^{(j)} \ne 0$ for $j =1, \dots, g$ (we
need this assumption in  Proposition \ref{perturbation} and
\ref{perturbation2}). The choice of the notation $\tr_j = (t_j,
\dots, t_g)$ in the deformation $P_{\tr_j}$ is related to the fact
that the terms $M_l (r,s)$ appearing in the expansion of
$P_{\tr_j}$ are expressed in terms of the monomial $x_l^r y_l^s$
at the local coordinates of the level $l$ of the embedded
resolution of $(C,0)$, for $j \leq l \leq g$.
 Notice that the polynomial $P_{\tr_1}$
determines any of the terms $P_{\tr_j}$ for $1 < j \leq g$, by
substituting
 $t_{1} = \dots = t_{j-1}
=0$ in $P_{\tr_1}$. Occasionally, we abuse of notation by denoting
$F$ by $P_{\tr_{g+1}}$.

\begin{Defi}
A {\em multi-semi-quasi-homogeneous  (msqh)   deformation} of the
plane branch $(C,0)$ is a family $C_{\tr_1} $ defined by,
$P_{\tr_1} =0 $ where $P_{\tr_1}$  is of the form (\ref{pg}), in a
Milnor ball of $(C,0)$. We say that $C_{\tr_1} $ is a
msqh-smoothing of $(C,0)$ if the curve $C_{\tr_1} $ is smooth and
transversal to the boundary of a Milnor ball for  $0 < t_1 \ll
\cdots \ll t_g \ll 1$.
\end{Defi}

\begin{Not}
We denote by $C_{\tr_l}$, or by $C_{\tr_l}^{(1)}$, the deformation
of $(C,0)$ defined by $P_{\tr_l}$ in a Milnor ball of $(C,0)$, for
$0 < t_l \ll \cdots \ll t_g \ll 1$ and $l=1, \dots, g$.

We denote by $C_{\tr_l}^{(j)} \subset Z_j$ the strict transform of
$C_{\tr_l}$ by the composition of toric maps $\Pi_{j-1} \circ
\cdots \circ \Pi_1$ and by $P_{\tr_l}^{(j)} (x_j, y_j)$ (resp. by
$\bar{P}_{\tr_l}^{(j)} (x_j, u_j) $) the polynomial defining
$C_{\tr_l}^{(j)}$ in the coordinates $(x_j, y_j)$ (resp. $(x_j,
u_j)$), for $2 \leq j \leq l \leq g$.
\end{Not}
These notations are analogous to those used for $C$ in Section
\ref{toric-res}, see  (\ref{strict-j}). In particular,  we have
that the result of substituting in $P_{\tr_j}^{(j)} (x_j, y_j)$,
the term $y_j$, by using Formula (\ref{xi-rel}) at level $j$, is
$\bar{P}_{\tr_l}^{(j)} (x_j, u_j)$.

\begin{Pro} \label{perturbation} {\rm (\cite{roots})}
\begin{enumerate}
\item[(i)] If $1 \leq j < l \leq g$ the curves $C_{\tr_l}^{(j)}$
and $C^{(j)}$ meet the exceptional divisor of $\Pi_{j-1} \circ
\cdots \circ \Pi_1$ only at the point $o_j \in E_j$ and with the
same intersection multiplicity $e_{j-1}$.

\item[(ii)] If $ 1 < j \leq g$ the curves $C_{\tr_j}^{(j)}$ meet
the exceptional divisor of $\Pi_{j-1} \circ \cdots \circ \Pi_1$
only at $e_j$ points of $E_j$, counted with multiplicity.
\end{enumerate}
\end{Pro}

\begin{Rem} \label{perturbation-np} {\rm (\cite{roots})}
\begin{enumerate}
\item[(i)] If $1 \leq j <  l \leq g$ the local Newton polygon of
$P_{ \tr_l}^{(j)} (x_j, y_j)$ and of $F^{(j)} (x_j, y_j)$
coincide.

\item[(ii)] If $1 <  j \leq g$ then $\Xi_j$ is a face of the
Newton polygons of $\bar{P}_{\tr_j}^{(j)} (x_j, u_j)$ and of
$P_{\tr_j}^{(j)} (x_j, y_j)$ corresponding to those terms which
are not divisible by $x_j$ (see Notations \ref{notation-res}).
\end{enumerate}
\end{Rem}

\begin{Pro} \label{perturbation2} (see \cite{roots})
If $1 < j \leq g$ then  we have that:
\begin{enumerate}
\item[(i)]The symbolic restriction of $P^{(j)}_{\tr_{j+1}} (x_{j},
y_{j}) $ to the edge $\Gamma_j$ of its local Newton polygon is of
the form:
\[
\a_j \prod_{s=1}^{e_{j}} \left( y_{j} ^{n_j} - (1 + \g_s^{(j)}
t_{j+1}^{e_{j+1} m_{j+1} } ) x_{j}^{m_j} \right),
\]
where $\a_j, \g_s^{(j)} \in \C \setminus \{ 0 \}$, for $s=1,
\dots, e_{j}$.

 \item[(ii)] The points of intersection of  $E_{j+1}$ with
$C_{\tr_{j+1}}^{(j+1)}$ are those with coordinates $x_{j+1} =0$
and
\begin{equation} \label{intersection2}
u_{j+1} = (1 + \g_s^{(j)} t_{j+1}^{e_{j+1} m_{j+1}})^{-1} , \mbox{
for } s=1, \dots, e_{j}.
\end{equation}
\end{enumerate}
\end{Pro}

\begin{Rem} \label{positive}
It follows from Proposition \ref{perturbation2} that those
peripheral roots of $P^{(j)}_{\tr_{j+1}} (x_{j}, y_{j})$, which
 are real,
are also positive for $0< t_{j+1}  \ll 1$.
\end{Rem}

When we say that $C_{\tr_l}^{(j)}$ defines a {\em deformation with
parameter $t_l$}, we mean for  $t_{l+1}, \dots t_{g}$ fixed.
Proposition \ref{msqh} motivates our choice of terminology in this
section.

\begin{Pro} \label{msqh} $C_{\tr_{j}}^{(j)}$ is a sqh-deformation with parameter $t_j$
of the singularity $(C_{\tr_{j+1}}^{(j)}, o_j)$ for $1 \leq j \leq
g$.
\end{Pro}
{\em Proof.} By Proposition \ref{perturbation} the curves $C_{
\tr_{j+1}}^{(j)}$ and $C_{\tr_j}^{(j)}$ intersect only the
irreducible component $E_j$ of the exceptional divisor of
$\Pi_{j-1} \circ \dots \circ \Pi_1$. By the construction of the
toric resolution this intersection is contained in the chart
$\C^2(\s_{j-1}) \subset Z_j$. By Lemma \ref{comb2} and the
definitions in Formula (\ref{pg}), if $C_{\tr_{j+1}}^{(j)}$ is
defined on the chart $\C^2(\s_{j-1})$ by $P_{\tr_{j+1}}^{(j)}
(x_j, y_j)$ then $C_{\tr_{j}}^{(j)}$ is defined by:
\[
P_{\tr_j}^{(j)}(x_j, y_j)  = \sum_{(r,s) \in \D_j^- \cap \Z^2}
A_{r,s}^{(j)} \, t_{j}^{\w_{j} (r,s)} \underline{u}
^{\underline{k}(r,s)} \, x_j^r \, y_j^s + P_{\tr_{j+1}}^{(j)}(x_j,
y_j),
\]
where for each $(r,s)$ the term $ \underline{u}
^{\underline{k}(r,s)}$ denotes the term $u_1^{k_1} \cdots
u_j^{k_j}$ of (\ref{ler}).  The elements $u_1, \dots, u_j$,
expanded in terms of  $x_j, y_j$, have constant term equal to one
by (\ref{units}). It follows from this that the local Newton
polygon of $P_{\tr_{j}}^{(j)}(x_j, y_j) $, with respect to the
variables $x_j$, $y_j$ and $t_j$, has only one compact face of
dimension two, which is equal to the graph of $\w_j$ on the Newton
diagram $\D_j$.
 $\Box$

Notice that the polynomial, $(P_{\tr_j}^{(j)})^{\hat{\D}_j}_{t_j
=1}$, defining the chart of the sqh-smoothing in Proposition
\ref{msqh} does not depend on $t_j$ (see Notation
\ref{notation-viro}).
We have that:
\begin{equation} \label{tri}
(P_{\tr_j}^{(j)})^{\hat{\D}_j}_{t_j =1} = \sum_{(r,s) \in \D_j^-
\cap \Z^2} A_{r,s}^{(j)} \, x_j^r \, y_j^s  +
(P_{\tr_{j+1}}^{(j)})^{\Gamma_j},
\end{equation}
where $(P_{\tr_{j+1}}^{(j)})^{\Gamma_j}$ is described by
Proposition \ref{perturbation2}.
\begin{Defi} \label{def-nd}
The msqh-deformation  $C_{\tr_1}$ is {\em real non degenerated} if
the polynomials $(P_{\tr_j}^{(j)})^{\hat{\D}_j}_{t_j =1}$ in
(\ref{tri}) are real non degenerated with respect to the polygon
$\D_j$, for $j=0, \dots, g-1$ (see Notations \ref{notation-viro}).
\end{Defi}

\begin{Pro} \label{msqh2} If the msqh-deformation $C_{\tr_1}$ is real non
degenerated then $C_{\tr_j}^{(j)}$ is a msqh-smoothing of the
singularity $(C^{(j)}, o_j)$. In particular, $C_{\tr_1}$ is a
msqh-smoothing of $(C,0)$.
\end{Pro}
{\em Proof.} We prove the Proposition by induction on $g$. If $g
=1$ then the assertion is a consequence of Definition \ref{sqh}.
Suppose $g>1$, then by the induction hypothesis $C_{\tr_2}^{(2)}$
is a msqh-smoothing of $(C^{(2)}, o_2)$. By Proposition
\ref{perturbation2} the polynomial $P_{\tr_2} (x, y)$, defining
the curve $C_{\tr_2}$, is non degenerated with respect to its
local Newton polygon. By Definition \ref{sqh} the deformation
$C_{\tr_1}$ is a sqh-smoothing of $(C_{ \tr_2}, 0)$ with parameter
$t_1$. It follows that $C_{\tr_1}$ defines then a msqh-smoothing
of $(C, 0)$ for $0 < t_1 \ll \cdots \ll t_g  \ll 1$. $\Box$

\subsection{Gluing the charts of msqh-smoothings} \label{glue}

In this Section we describe the patchwork  of the charts of
$(P_{\tr_j}^{(j)})^{\hat{\D}_j}_{t_j =1}$ and of
$P_{\tr_{j+1}}^{(j)}$ at the level $j$ of the toric resolution,
under some geometrical assumptions. We begin by some Lemmas for
sqh-smoothings.

\begin{Lem}  \label{compatible-oscillation}
Let us consider a semi-quasi-homogeneous smoothing $C_t$ defined
by $P_t (x, y) =0$ of $(C, 0)$ (see Notations \ref{notation-sqh}).
Set $\Lambda_0$ for the Newton polygon of $P_0$. Suppose that the
following statements hold:
\begin{enumerate}
 \item[(i)] The
chart $\mbox{\rm Ch} _{\Lambda_0/\sim } (P_0)$ is in maximal
position with respect to $\tilde{\Gamma}/\sim$.

\item[(ii)]  The chart $\mbox{\rm Ch} _{\D/ \sim} (
P^{\hat{\D}}_{{1}} )$ is in maximal position with respect to
$\tilde{\Gamma}/\sim$.

\item[(iii)] The order of the peripheral roots of $P^{\Gamma}$
coincide when viewed in suitable arcs of the charts  $\mbox{\rm
Ch} _{\Lambda_0/\sim } (P_0)$ and $\mbox{\rm Ch} _{\D/ \sim} (
P^{\hat{\D}}_{{1}} )$.
\end{enumerate}
We label the   peripheral roots $\a_1, \dots, \a_e \in
\tilde{\Gamma}/\sim$  of $P^{\Gamma}$ with the order induced by
the arcs the charts in {\rm (i)} and  {\rm (ii)}. We denote by
$\mathbf{a}_{{{k}}} {(j, j+1)}$ the arc of the chart in (k), for
$k ={{i,ii}}$, which joins the peripheral roots $\a_j$ and
$\a_{j+1}$, for $j =1, \dots, e-1$. Then exactly one of this two
statements is verified:
\begin{enumerate}
 \item[(a)]
 In the patchwork of the
charts of $P^{\hat{\D}}_{{1}} ( (x,y)) $ and of $P_0$ the arcs
$\mathbf{a}_{{{i}}} {(j, j+1)}$ and $\mathbf{a}_{{{ii}}} {(j,
j+1)}$ glue into an oval intersecting $\tilde{\Gamma}$, for $j =1,
\dots, e-1$.

\item[(b)] In the gluing of the charts of $ P^{\hat{\D}}_{{1}}
(\r_\Gamma (x,y)) $ and of $P_{0}$ the arcs $\mathbf{a}_{{{i}}}
{(j, j+1)}$ and $\r_{\Gamma} (\mathbf{a}_{{{ii}}} {(j, j+1)})$
glue into an oval intersecting $\tilde{\Gamma}$, for $j =1, \dots,
e-1$.
\end{enumerate}
\end{Lem}
{\em Proof.} If $e =1$ there is nothing to prove, hence we suppose
that $e>1$. By forgetting the relation $\sim$, we obtain two
symmetric copies of $\a_1, \dots, \a_e$ in $\tilde{\Gamma}$, by
the action of the symmetry $\r_\Gamma$ (see notations of Section
\ref{toric}). We denote by $\mathbf{a}_{k}{(1)}$  (resp. by
$\mathbf{a}_{{{k}}}{(e)}$) the arc of chart in (k), for $k
=\mbox{{\rm i,ii}}$, which intersects $\mathbf{a}_{{{k}}}{(1, 2)}$
at $\a_1$ (resp. $\mathbf{a}_{{{k}}}{(e-1, e)}$ at $\a_e$).

By Theorem \ref{vlocal} if the arcs $\mathbf{a}_{{{i}}}{(1, 2)}$
and $\mathbf{a}_{{{ii}}}{(1, 2)}$ glue up in an oval intersecting
$\tilde{\Gamma}$  in the gluing of the charts of
$P^{\hat{\D}}_{{1}}  (x,y) $ and of $P_{\Lambda_0}$, then the
assertion of the Lemma holds for $\r = Id$, otherwise the
assertion of the Lemma holds for $\r = \r_\Gamma$.
 $\Box$

\begin{figure}[htbp]
$$\epsfig{file=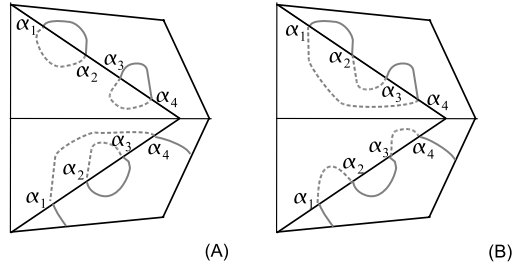, height= 35 mm}$$
\caption{A represents regular intersection of charts
\label{maximal-pos}}
\end{figure}

 The dotted style curve in  Figure \ref{maximal-pos} represents  the two
possibilities for the chart of $P^{\hat{\D}}_{{1}}$ with good
oscillation. Cases (A) and (B) correspond to assertion (a) and (b)
respectively, where in this case  $\r_{\Gamma}$ is the symmetry
$\r_{\Gamma} (r, s) = (r, -s)$.

 The following terminology is introduced,
 with a slightly different meaning, in \cite{KRS}:

\begin{Defi} \label{maximal}
We say that the charts of $P^{\hat{\D}}_{{1}} $ and of $P_0$,
associated with the sqh-smoothing $C_t$, have regular intersection
along $\Gamma$ if the statement (a) of Lemma
\ref{compatible-oscillation}  hold.
\end{Defi}

\begin{Lem} \label{compatible-oscillation2}
Suppose that in the Lemma \ref{compatible-oscillation} statement
(a) holds. If in addition there exists (resp.  there does not
exist) a connected component of the chart $\mbox{\rm
Ch}^*_{\Lambda_0} (P_0)$   bounded by $\a_1$ and $\a_e$ then in
the gluing of the charts of $P^{\hat{\D}}_{{1}}  ( (x,y)) $ and of
$P_{0}$ the arcs $\mathbf{a}_{{{i}}}{(1)}$,
$\mathbf{a}_{{{ii}}}{(1)}$, $\mathbf{a}_{{{i}}}{(e)}$ and
$\mathbf{a}_{{{ii}}}{(e)}$ are contained  (resp. are not
contained) in an oval
 intersecting $\tilde{\Gamma}$.
\end{Lem}
{\em Proof.} Notice that by construction the arcs
$\mathbf{a}_{{{i}}}{(1)}$ glue with $\mathbf{a}_{{{ii}}}{(1)}$
(resp. for $\mathbf{a}_{{{i}}}{(e)}$ and
$\mathbf{a}_{{{ii}}}{(e)}$).

The arcs $\mathbf{a}_{{{ii}}}{(1)}$ and $\mathbf{a}_{{{ii}}}{(e)}$
are in the same connected component of the chart $\mbox{\rm Ch}
_{\D/ \sim} (  P^{\hat{\D}}_{{1}} )$ since  this chart is compact
and hence each connected component is an oval in particular the
one which is in maximal position with respect to
$\tilde{\Gamma}/\sim$.

The statement follows easily from these observations and the
hypothesis. $\Box$

\begin{figure}[htbp]
$$\epsfig{file=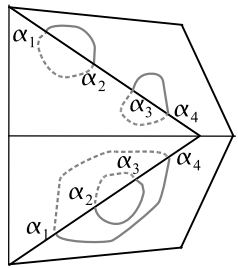, height= 35 mm}$$
\caption{Regular intersection of charts\label{maximal-pos2}}
\end{figure}

Figures  \ref{maximal-pos} case (A) and \ref{maximal-pos2}
represent the two possibilities indicated in Lemma
\ref{compatible-oscillation2} when the chart $P^{\hat{\D}}_{{1}}$
has good oscillation with respect to $\tilde{\Gamma}/\sim$.

Denote by ${\Lambda}_{2}$ the Newton polygon of  $P_{\tr_{2}}
(x_1, y_1)$. Notice that ${\Lambda}_{2}$ contains the Newton
polygon of $F^{(1)}(x_1, y_1)$ and shares with it the common  edge
$\Gamma_1$ by Remark \ref{perturbation-np} (i). We denote by
$\bar{\Lambda}_{2} ^{(2)}$ the Newton polygon of
$\bar{P}_{\tr_{2}} ^{(2)} (x_{2}, u_{2})$.

\begin{Pro} \label{inter}
Let $C_{\tr_1}$ be a real non degenerated msqh-smoothing of
$(C,0)$. If the smoothing $C_{\tr_{2}}^{(2)}$  of $(C^{({2})},
o_{2})$  is in maximal position (resp. has good oscillation) with
respect to the line $E_{2}$ then in a neighborhood of
$\tilde{\Gamma}_{{1}}/\sim$,
 the chart $\mbox{\rm
Ch}_{\Lambda_{2}
/\sim} (P_{\tr_{2}}
(x_{1}, y_{1}))$ is in maximal position (resp. has good
oscillation) with respect to $\tilde{\Gamma}_{{1}}/\sim$.
\end{Pro}
{\em Proof.} By (\ref{strict-j}) we have that:
\[
P_{\tr_{2}}
\circ \p_{{1}} = \mbox{{\rm Exc}}   (
P_{\tr_{2}}
, \p_{1} ) \cdot  \bar{P}_{\tr_{2}} ^{(2)}
(x_{2}, u_{2}),
\]
where $\mbox{{\rm Exc}}  ( P_{\tr_{2}} , \p_{1} ) $ is a monomial
in $x_{2}$ and  $u_{2}$.  Notice that $\bar{P}_{\tr_{2}} ^{(2)}
(0, u_{2})$ is a polynomial of degree $e_{1}$ with non zero
constant term by Proposition \ref{perturbation} (ii). Denote by
$\phi$ the composition of the $\mbox{{\rm SL}} (2, \Z)$
transformation corresponding to (\ref{mm-i}) with the translation
induced by  the exponent of the monomial $\mbox{{\rm Exc}}  (
P_{\tr_{2} } , \p_{1} ) $. Then we have that $\phi$ is an
isomorphism of triples
\begin{equation} \label{triple}
\left( \bar{\Lambda}_{2} ^{(2)}, \Xi_{2}, \mbox{\rm Ch} _{
\bar{\Lambda}_{2} ^{(2)} /\sim } (\bar{P}_{\tr_{2}} ^{(2)} (x_{2},
u_{2}) ) \right) \stackrel{\phi}{\longrightarrow} \left(
{\Lambda}_{2} , \Gamma_{{1}}, \mbox{\rm Ch} _{ {\Lambda}_{2}
 /\sim } (P_{\tr_{2}}
 (x_{{1}},
y_{{1}})) \right).
\end{equation}
In other terms $\p_{1}$ induces an isomorphism of the toric
surfaces associated to these Newton polygons which maps the one
dimensional orbit associated to $\Gamma_{{1}}$ to the orbit
associated to $\Xi_{2}$. More precisely, a point of the orbit of
$\Gamma_{{1}}$ is defined by the vanishing of $y_{{1}}^{n_{1}} -
\theta x_{{1}}^{m_{1}}$ for some $\theta \in \C^*$ and then by the
computations of Section \ref{toric-res} it corresponds to the
point of coordinates $x_{2} =0$ and $u_{2} = \theta^{-1}$,
identified with the orbit associated to $\Xi_{2}$. This implies
that the following assertions are equivalent:
\begin{enumerate}
\item[(a)] In a neighborhood of the line $\tilde{\Xi}_{2} /\sim$,
the chart $\mbox{\rm Ch} _{ \bar{\Lambda}_{2} ^{(2)} /\sim } (
\bar{P}_{\tr_{2}} ^{(2)} (x_{2}, u_{2})) $ is in maximal position
(resp. has good oscillation) with respect to the line
$\tilde{\Xi}_{2} /\sim$.

\item[(b)] In a neighborhood of the line,
$\tilde{\Gamma}_{{1}}/\sim$,  the chart $\mbox{\rm Ch} _{
{\Lambda}_{2} /\sim } (P_{\tr_{2}} (x_{{1}}, y_{{1}})) $ is in
maximal position (resp. has good oscillation) with respect to
$\tilde{\Gamma}_{{1}}/\sim$.
\end{enumerate}

If the smoothing $C_{\tr_{2}}^{(2)}$ is in maximal position (resp.
has good oscillation) with respect to the line $E_{2}$ then there
is an arc  of the smoothing  $\R C_{ \tr_{2}}^{(2)}$ which
contains  the points $E_{2} \cap C_{\tr_{2}}^{(2)}$ and satisfies
the geometrical hypothesis with respect to the line $E_{2}$. By
Proposition \ref{perturbation} this arc does not intersect the
coordinate line corresponding to $u_{2} =0$ in the chart. This
implies that (a) holds and hence (b) holds. $\Box$

\begin{Lem} \label{em}
Let $C_{\tr_1}$ be a real non degenerated msqh-smoothing of
$(C,0)$. Suppose that the charts $\mbox{\rm Ch}^*_{\D_{1}}
(P^{\hat{\D}_1}_{t_1 =1} (x_{1}, y_{1}))$ and $\mbox{\rm
Ch}^*_{\Lambda_{2} } (P_{\tr_{2}} (x_{1}, y_{1}))$ have regular
intersection along $\Gamma_{1}$ (see Definition \ref{maximal}). If
the connected component of $C_{\tr_{2}}^{(2)}$  which meets
$E_{2}$ is an oval (resp. it is not an oval),  then there are
precisely $e_{1}$ ovals  (resp. $e_{1} -1$ ovals and one non
closed component) which intersect $\tilde{\Gamma}_{1}$ in the
Patchwork of the charts $\mbox{\rm Ch}^*_{\Lambda_{2} }
(P_{\tr_{2}}^{({1})})$ and of $\mbox{\rm Ch}^*_{\D_{1}}
(P^{\hat{\D}_1}_{t_1 =1})$, which describes the smoothing
$C_{\tr_{1}} $ of $C_{\tr_{2}} $. In addition, if the chart
$\mbox{\rm Ch} _{ {\Lambda}_{2} /\sim } (P_{\tr_{2}} (x_{{1}},
y_{{1}})) $ has good oscillation) with respect to
$\tilde{\Gamma}_{{1}}/\sim$ and if $P^{\hat{\D}_1}_{t_1 =1}$
defines a Harnack curve in $Z(\D_{1})$, then the signed
topological type of the chart of $P^{\hat{\D}_1}_{t_1 =1}$ is
unique.
\end{Lem}
{\em Proof.} Notice that by hypothesis and Proposition \ref{inter}
the smoothing $C_{\tr_{2}}^{(2)}$ is in maximal position with
respect to the line $E_{2}$. We have that the chart of $\mbox{\rm
Ch}_{\Lambda_{2} /\sim} (P_{\tr_{2}} ) $
 is in maximal position with respect to the line $\Gamma_{1} /\sim$.
We denote by $\a_1, \dots \a_{e_{1}}$ the peripheral roots of
$P_{\tr_{2}} ^{\Gamma_{1}}$, labelled with respect to the order in
the charts $\mbox{\rm Ch}_{\Lambda_{2} /\sim} (P_{\tr_{2}} )$ and
$\mbox{\rm Ch}_{\D_{1}/\sim} (P^{\hat{\D}_1}_{t_1 =1})$. Then
there exist $e_{1} -1 $ ovals meeting $\tilde{\Gamma}_{1}$ in the
patchwork of these charts by Lemma \ref{compatible-oscillation}.
In addition, if the connected component of  $\R C_{\tr_{2}}^{(2)}$
which meets $E_{2}$ is an oval (resp. is not an oval), then the
hypothesis of Lemma \ref{compatible-oscillation2} are satisfied
and the conclusion follows. For the second statement by
Proposition \ref{signed-topological}  we have that there are two
possible signed topological types for $P_{\D}$ with prescribed
symbolic restriction to the face $\Gamma_1$. By Lemma
\ref{compatible-oscillation} only one of these two types induces
regular intersection along the edge $\Gamma_1$. $\Box$

\begin{figure}[htbp]
$$\epsfig{file=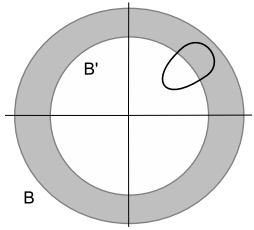, height= 30 mm}$$
\caption{A mixed oval \label{mixed-oval}}
\end{figure}

We call the ovals described by Lemma \ref{em}  {\em mixed ovals of
depth $1$}. We call {\em ovals of depth $1$},  those which appear
in the chart of $P^{\hat{\D}_1}_{t_1 =1} $ (see (\ref{tri}))
 but do not cut $\tilde{\Gamma}_1$. In Figure \ref{mixed-oval} we represent
a mixed oval of depth one; the ball $B$ is a Milnor ball for
$(C_{\tr_{2}},0)$,  the segment of the oval in small ball $B'$
corresponds to an arc of the chart of $P^{\hat{\D}_1}_{t_1 =1} $,
while the segment of the oval in $B \setminus B'$ corresponds to
an arc of the chart $\mbox{\rm Ch}_{\Lambda_{2} /\sim}
(P_{\tr_{2}} ) $.

\begin{Defi} \label{depth} Let  $C_{\tr_1}$ be a msqh-smoothing  of
$(C,0)$. An oval $O$ of  $C_{\tr_1}$  is of depth $j$ (resp. a
mixed oval of depth $j$) if there exists an oval $O_j$ of depth
$1$ (resp. a mixed oval of depth $1$) of the smoothing
$\C_{\tr_j}^{(j)}$ of $(C^{(j)}, o_j) $ such that $E_j \cap O_j =
\emptyset $ and such that $O$ arises as a slight perturbation of
the oval $\Pi_1 \circ \cdots \circ \Pi_{j-1} (O_j)$ of
$C_{\tr_j}$, for $1 \leq j \leq g$.
\end{Defi}

\subsection{Maximal, Harnack and multi-Harnack smoothings}

We introduce the following notions for a real non degenerated
msqh-smoothing    $C_{\tr_1} $  of a real plane branch $(C,0)$. By
Proposition \ref{msqh2} if $C_{\tr_1}$ is a non degenerated
msqh-smoothing of $(C,0)$ then $C_{\tr_j}^{(j)}$ is also a non
degenerated msqh-smoothing of $(C^{(j)}, o_j)$, for $1 \leq j \leq
g$.
\begin{Defi} \label{above}
\begin{enumerate}
\item[(i)] $C_{\tr_1} $   is a $M$-{\em smoothing} if the number
of ovals  in a Milnor ball of $(C,0)$ is equal to $\frac{1}{2} \mu
(C)_0$.

\item[(ii)] A $M$-smoothing  $C_{\tr_1} $ is {\em Harnack} if it
has good oscillation with respect to the coordinate axis.
\item[(iii)] A $M$-smoothing $C_{\tr_1}$ is {\em multi-Harnack} if
$C_{\tr_j}^{(j)}$ is a Harnack $M$-msqh-smoothing of $(C^{(j)},
o_j)$, for $1 \leq j \leq g$.
\end{enumerate}
\end{Defi}
In Definition \ref{above} (iii) the Harnack condition is
considered with respect  to the coordinate lines defined by the
coordinates $(x_j, y_j)$, see Section \ref{toric-res} for notations.

The following result describes inductively $M$-msqh-smoothings and
Harnack msqh-smoothings of $(C,0)$ when $g
>1$.
\begin{The} \label{main} Let $C_{\tr_1}$  be a non
degenerated msqh-smoothing of the real plane branch $(C, 0)$. We
introduce the following conditions:
\begin{enumerate}

\item[(i)] $C_{\tr_2}^{(2)}$ is a $M$-msqh-smoothing of $(C^{(2)},
o_2)$ in maximal position with respect to the exceptional divisor
$E_2$.

\item[(ii)]  $C_{\tr_1}$ defines a $M$-sqh-smoothing of the
singularity $(C_{\tr_2}, 0)$ with parameter $t_1$.

\item[(iii)]  The charts of $P^{\hat{\D}_1}_{t_1 =1}$ and of
$P_{\tr_2}$ have regular intersection along $\Gamma_1$.
\end{enumerate}
Then, the deformation $C_{\tr_1}$ is a $M$ (resp Harnack)
msqh-smoothing of $(C, 0)$ if and only if conditions {\rm (i),
(ii)} and {\rm (iii)} hold (resp. and in addition $C_{\tr_1}$ has
good oscillation with respect to the coordinate axis).
\end{The}
{\em Proof.} We prove first that if conditions {\rm (i)} and {\rm
(ii)} hold then $C_{\tr_1}$ is a $M$-msqh-smoothing of $(C, 0)$.
The number of ovals of the msqh-smoothing $C_{\tr_2}^{(2)}$ is
equal to
\begin{equation} \label{arr}
\frac{1}{2} \mu (C^{(2)})_{o_2} \stackrel{(\ref{mu-delta-j})}{=}
 \sum_{j=2}^{g}  \left( \# ( \mbox{{\rm int}} \D_j \cap \Z^2 )
+ e_{j} -1 \right).
\end{equation}
Since $C_{\tr_2}^{(2)}$ is in maximal position with respect to
$E_2$ there is only one connected component $A$ of the smoothing
$C_{\tr_2}^{(2)}$ which intersect $E_2$ in $e_1$ different real
points. When we apply the toric morphism $\p_1$ we get a
deformation $(C_{\tr_2},0)$ of $(C, 0)$, defined by $P_{\tr_2} (x,
y)$. Notice that $(C_{\tr_2},0)$ is singular at $0$. By
Proposition \ref{perturbation2} the singularity $P_{\tr_2} (x,
y)=0$ is real non degenerated with respect its local Newton
polygon. The image $A'$  by $\p_1$ of the component $A$ passes
through the origin and is the only connected component of $\R
C_{\tr_2}$ with this property. The deformation  $C_{\tr_2}$  has
the same number of ovals as $C_{\tr_2}^{(1)}$, which are of depth
$> 1$.

 By hypothesis (ii) we
have that $C_{\tr_1}$ is a $M$-sqh-smoothing of $C_{\tr_2}$ with
parameter $t_1$ hence it yields $\# ( \mbox{{\rm int}} \D_{1} \cap
\Z^2 )$ ovals of depth $1$. By Proposition \ref{inter} the chart
of $P_{\tr_2} (x, y)$ with respect to its Newton polygon
$\Lambda_{2}$, is in maximal position with respect to
$\Gamma_1/\sim$. By hypothesis (iii) we are in the situation
described by Lemma \ref{em}: if $A$ is an oval (resp. is not) the
image of $A$ in the chart of $P_{\tr_2} (x, y)$ patchwork with the
chart $\mbox{\rm Ch}^*_{\D_{1}} (P^{\hat{\D}_1}_{t_1 =1})$
providing $e_1$ mixed ovals of depth $1$ (resp.  $e_1 -1$ mixed
ovals of depth $1$). It follows that the msqh-smoothing
$C_{\tr_1}$ has the maximal number of   ovals (see
(\ref{mu-delta-j})).

Conversely, suppose that $C_{\tr_1}$ is a msqh-smoothing. Since
$C_{\tr_1}$ defines a sqh-smoothing of $(C_{\tr_2}, 0)$ with
parameter $t_1$,  there are at most $\# ( \mbox{{\rm int}} \D \cap
\Z^2 )$ ovals of depth $1$.

 If there are $r \geq 1 $ components of the smoothing $C_{\tr_2}^{(2)}$
 which intersect $E_2$ then
we prove that the maximal number of ovals of the smoothing $C_{
\tr_1}$
is bounded below by:
\[
 -r + 1 + \sum_{j=0}^{g-1}  \left( \# ( \mbox{{\rm int}} \D_j
\cap \Z^2 ) + e_{j+1} -1 \right)  .
\]
If no component of the smoothing  $C_{ \tr_2}^{(2)}$ intersects
$E_2$ then there is no mixed oval of depth $1$  for the smoothing
$C_{\tr_2}$ and therefore  $C_{ \tr_1}$ is not a
$M$-msqh-smoothing by (\ref{mu-delta-j}). Therefore there exist
$r\geq 1$ components of the smoothing $C_{\tr_2}^{(2)}$, each one
cutting $E_2$ in $s_r \geq 1$ real points. We argue as in Lemma
\ref{em}: if such a component
 is not an oval
then it leads to $s_r -1$ mixed ovals;
 if this component is an oval then it contributes with
$s_r$   mixed ovals, but we also loose one oval of depth $2$. It
follows that  if $C_{\tr_1}$ is a $M$-msqh-smoothing then $r=1$
and $s_1 = e_1$, i.e., $C_{\tr_2}^{(2)}$ is in maximal position
with respect to the line $E_2$ and assertion (iii) hold. $\Box$

\begin{Cor}      \label{multi-h2}
Let $C_{\tr_1}$  be a real non degenerated msqh-smoothing of the
real plane branch $(C, 0)$. We introduce the following conditions:
\begin{enumerate}

\item[(i)] $C_{\tr_2}^{(2)}$ is a multi-Harnack smoothing of
$(C^{(2)}, o_2)$.

\item[(ii)]  $C_{\tr_1}$ defines $M$-sqh-smoothing of the
singularity $(C_{\tr_2}, 0)$ with parameter $t_1$, in maximal
position with respect to the coordinate axis.

\item[(iii)] The charts of $P^{\hat{\D}_1}_{t_1 =1}$ and of
$P_{\tr_2}$ have regular intersection along $\Gamma_1$.
\end{enumerate}
Then the deformation $C_{\tr_1}$ is a multi-Harnack smoothing of
$(C, 0)$ if and only if conditions {\rm (i), (ii)} and {\rm (iii)}
hold.
\end{Cor}
{\em Proof.} It follows by induction on $g$ from Theorem
\ref{main} and Theorem \ref{1par}. $\Box$

\subsection{Existence of multi-Harnack smoothings}

In this section we prove the existence of a multi-Harnack
smoothing of a real plane branch $(C,0)$. We will use the
following observation of Viro.

\begin{Rem} \label{vsqh} {\rm (see \cite{Vpw} page 19)} If $(C,0)$ is
defined by a non degenerated semi-quasi-homogeneous polynomial
then any smoothing of $(C,0)$  constructed by patchworking is {\em
topologically equivalent}, in a stratified sense with respect to
the boundary of a Milnor ball and the coordinate axis, to a
sqh-smoothing.
\end{Rem}

\begin{The}
Any real plane branch $(C, 0)$ has a multi-Harnack smoothing.
\end{The}
{\em Proof.} We construct a multi-Harnack smoothing $C_{\tr_1}$ of
$(C, 0)$ by induction on $g$. If $g=1$ by Proposition
\ref{compatible} there exists a Harnack smoothing of $(C, 0)$ and
by Remark \ref{vsqh} we construct from this a topologically
equivalent sqh-smoothing, which is also Harnack. Suppose the
result true for $g-1$. Then by induction hypothesis we have
constructed a multi-Harnack smoothing $C_{\tr_2}^{(2)} $ of
$(C^{(2)}, o_2)$. By Proposition \ref{msqh} the deformation
$C_{\tr_2}$ is defined by a sqh-polynomial $P_{\tr_2} (x, y)$ with
peripheral roots of the same sign by Remark \ref{positive}. Then
we can apply Proposition \ref{compatible} and Remark \ref{vsqh} to
construct Harnack sqh-smoothings $C_{\tr_1}$ of the singularity
$(C_{\tr_2}, 0)$ with parameter $t_1$, such that the charts of
$P^{\hat{\D}_1}_{t_1 =1}$  and of $P_{\tr_2}$  have regular
intersection along $\Gamma_1$, by Lemma
\ref{compatible-oscillation}. The result follows by Corollary
\ref{multi-h2}. $\Box$

\subsection{The topological type of a multi-Harnack smoothing}

The definition of the topological type (resp. signed topological
type) of a msqh-smoothing $C_{\tr_1}$ of a real plane branch is
the same as in the  $1$-parametrical case (see Definition
\ref{topological2}).

\begin{The} \label{main-top}
Let $(C, 0)$ be a real branch. The topological type of
multi-Harnack smoothings of $(C, 0)$ is unique. There is at most
two signed topological types of multi-Harnack smoothings of $(C,
0)$. These types depend only on the complex equisingularity class
of $(C, 0)$.
\end{The}
{\em Proof.}   We prove the result by induction on $g$. For $g=1$
there is a unique topological type of Harnack smoothing  and two
signed topological types   by Theorem \ref{1par}. Suppose $g
>1$. We consider a multi-Harnack smoothing $C_{\tr_1}$ of $(C,0)$.
It is easy to see by Corollary \ref{multi-h2} and induction that
$C_{\tr_{g}}^{(g)}$ defines a sqh-Harnack smoothing of $(C ^{(g)},
o_{g})$. By construction the singularity $(C ^{(g)}, o_{g})$ is
non degenerated with respect to its Newton polygon in the
coordinates $(x_{g}, y_{g})$. Since the smoothing
$C_{\tr_{g}}^{(g)}$ is Harnack then there is only one topological
type and two signed topological types  (see Proposition
\ref{signed-topological}). We prove that there are at most two
signed topological types (resp. exactly one topological type) of
multi-Harnack smoothings, each one determined by the
equisingularity class of $(C,0)$ and an initial choice for the
sign type of the Harnack smoothing $C_{\tr_g}^{(g)}$: to be
normalized or not.

By induction on $g$ we  assume that the topological type (resp.
the signed topological type)  of the msqh-smoothing $C_{
\tr_2}^{(2)}$ of $(C^{(2)}, o_2)$ depends only on the
characteristic pairs of $(C^{(2)},o_2)$ with respect to the
coordinates $(x_2, y_2)$. Hence  it is determined by the complex
equisingularity class of $(C, 0)$. Since $\Pi_1$ is an isomorphism
over $\C^2 \setminus \{0 \}$  we deduce that the topological type
(resp. the signed topological type) of the msqh-deformation
$C_{\tr_2}$ is determined inside a Milnor ball $B$ of $(C,0)$, for
$0< t_2 \ll \cdots \ll t_g \ll 1$. Notice that $C_{\tr_2}$  is a
singular curve with $e_1$ branches at the origin. Let $B' \subset
{B}$ be a Milnor ball for the singularity $(C_{\tr_2}, 0)$. The
radius of $B'$  depends on $\tr_2$ and $B'$ is contained in the
interior of $B$. By construction the msqh-smoothing $C_{\tr_1}$ is
built by the patchwork of the charts of $P^{\hat{\D}_1}_{t_1 =1}$
and of $P_{\tr_2}$ described in Section \ref{glue}. By induction
hypothesis the topological type (resp. the signed topological
type)  of $C_{\tr_2}$ is fixed in $\R B \setminus \R {B}'$. We can
assume that the embedded topology of the chart
$\mbox{Ch}_{\Lambda_{2}} ((P_{\tr_2} )_{\Lambda_{2}})$ is
determined. By Corollary \ref{multi-h2} it is enough to prove that
there is a unique (resp. signed) topological way to patchwork the
chart of $P^{\hat{\D}}_{t_1 =1}$, in such a way that $C_{\tr_1}$
defines a multi-Harnack smoothing.

Corollary \ref{multi-h2} (ii) together with Theorem \ref{1par}
implies that $P^{\hat{\D}_1}_{t_1 =1}$ defines a Harnack curve in
the toric surface $Z(\D_{1})$. By Proposition
\ref{signed-topological} there are two possible signed topological
types for the chart of $P^{\hat{\D}_1}_{t_1 =1}$ (which are
related by the symmetry $\r_{\Gamma_1}$). By Corollary
\ref{multi-h2} the charts of $P_{\tr_2}$ and of
$P^{\hat{\D}_1}_{t_1 =1}$ have regular intersection along the edge
$\Gamma_1$. By Lemma \ref{compatible-oscillation2} this condition
holds for only one of the two possible signed types of Harnack
curves $P^{\hat{\D}_1}_{t_1 =1}$.

 The signed topological type of
the chart of $P^{\hat{\D}_1}_{t_1 =1}$ is uniquely determined and
it depends only on $\D_{1}$. By the induction hypothesis the
topological type of $C_{\tr_2}$ is fixed in $\R B \setminus \R
{B}'$ hence there is a unique topological type and at most two
signed topological of multi-Harnack smoothings of the real plane
branch $(C,0)$. The topological type only depend on the sequence
of triangles $\D_j$, hence on the equisingularity class of the
branch.  $\Box$

\subsection{Positions and scales of the ovals of multi-Harnack
smoothings}\label{taille}

We deduce some consequences on the position in the quadrants of
$\R^2$ and on the  scales of the ovals of a msqh-smoothing , in
particular in the multi-Harnack case.

\subsubsection{Multi-scaled structure of msqh-smoothings}

Let $C_{\tr_1}$ be a msqh-smoothing of a real plane branch
$(C,0)$. Recall that the notion of oval of depth $j$ of
$C_{\tr_1}$ was introduced in Definition \ref{depth}.

\begin{Pro} \label{size-j}
The size of an oval $O$ of depth $j$ in a msqh-smoothing
$C_{\tr_1}$ is equal to $(t_j^{e_1 n_1}, t_j^{e_1 m_1})$.
\end{Pro}
{\em Proof.} By Proposition \ref{size-sqh} the size of the oval
$O$ in the coordinates $(x_j, y_j)$ is equal to $(t_j^{e_{j} n_j}
,t_j^ {e_j m_j })$. The image of $O$ by $\p_{j-1}$, described by
(\ref{mm-i}), is an oval of $C_{\tr_j}^{(j-1)}$ of size $(t_{j}
^{e_j n_j n_{j-1}}, t_j^{e_j n_j m_{j-1}})$ since the function
$u_j \sim 1$ on the oval by (\ref{xi-rel}). By definition the oval
$\p_{j-1} (O)$ is not contained in the connected component of $\R
C_{\tr_j}^{(j-1)}$ which passes by $o_{j-1}$. By Proposition
\ref{size-sqh} and Theorem \ref{vlocal} the smoothing
$C_{\tr_{j-1}}^{(j-1)}$ of the singularity $(C_{\tr_j}^{(j-1)},
o_{j-1})$ only performs small perturbations of the ovals of
$C_{\tr_j}^{(j-1)}$ and hence the size of corresponding ovals with
respect to the coordinates $(x_{j-1}, y_{j-1})$ is the same.

It follows that $\p_{j-1} (O)$ is slightly deformed into an oval
$O_{j-1}$ of $C_{\tr_{j-1}}^{(j-1)}$ such that $E_{j-1} \cap
O_{j-1} = \emptyset $. By induction we find that the oval $O$ is
of size
\[
(t_j^{e_j n_j \cdots n_1}, t_j ^{e_j n_j \cdots n_2 m_1})  = ( t_j
^{e_1 n_1}, t_j^{e_1 m_1}). \quad \Box
\]

\begin{Pro} \label{non-nested}
If $C_{\tr_1}$ is a msqh-smoothing of the real plane branch
$(C,0)$ then any pair of ovals of depth $j, j'$ with $j < j'$ is
not nested.
\end{Pro}
{\em Proof.} By the Definition \ref{depth} it is enough to prove
it when $j =1$. Let $B$ denote a Milnor ball for the singularity
$(C, 0)$ and $B' \subset {B}$ a Milnor ball for the non
degenerated singularity $(C_{\tr_{2}}, 0)$. The ball $B'$ is
contained in the interior of $B$ and the radius of $B'$ depends on
the parameters. Notice that by Theorem \ref{vlocal} the smoothing
$C_{\tr_{1}}$ of the singularity $(C_{\tr_{2}}, 0)$ is constructed
by the patchwork of the charts of $P^{\hat{\D}_1}_{t_1 =1}$ and of
$P_{\tr_2}$. We have that topologically in a stratified sense with
respect to the coordinate axis and the boundary of $B'$, we
replace the pair $(\R {B}', C_{\tr_2} \cap \R {B}')$ by
$(\tilde{\D}, \mbox{Ch}_\D (P^{\hat{\D}}_{t_1 =1} )$, by
identifying $(\R {B}',
\partial \R   {B}'  )$   with $(\tilde{\D}, \partial
\tilde{\D})$, while in ${B} \setminus B'$ the curves $\R
C_{\tr_1}$ and $\R C_{\tr_{2}}$ remain isotopic. This implies that
in ${B} \setminus B'$ the curve $\R C_{\tr_1}$ contains arcs of
the mixed ovals of depth $1$ and ovals of depth $> 1$. $\Box$

\begin{Rem}
A msqh-smoothing of the real plane branch may have nested ovals
(see Example \ref{nested}).
\end{Rem}
\subsubsection{The multi-Harnack case}

\begin{figure}[htbp]
$$\epsfig{file=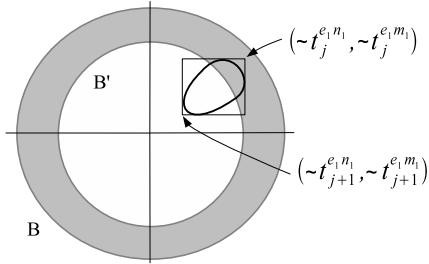, height= 35 mm}$$
\caption{Scale of a mixed oval of depht $j$ \label{mixed-oval2}}
\end{figure}

Let $C_{\tr_1}$ be a multi-Harnack smoothing of the real plane
branch $(C,0)$. By the proof of Theorem \ref{main-top}, the ovals
of $C_{\tr_j}^{(j)}$ do not cut $E_j$ for $j=1, \dots, g$. It
follows that each oval of $C_{\tr_j}^{(j)}$ is either an oval
(resp. a mixed oval) of depth $j$ for some $j$. There are
precisely $e_j -1$ mixed ovals of depth $j$ and $\# \mbox{\rm int}
(\D_j \cap \Z^2)$ ovals of depth $j$. The ovals of a multi-Harnack
smoothing are not nested by Proposition \ref{top-harnack}, Theorem
\ref{1par} and Proposition \ref{non-nested}. We complement the
information on the size of ovals of Section \ref{taille} by
indicating the form of the boxes which contain the mixed ovals of
depth $j$. The proof is analogous to that of Proposition
\ref{size-j} which indicates the size of ovals of depth $j$ (see
Figure \ref{mixed-oval2}).
\begin{Pro} \label{mixed size-j}
A mixed oval of depth $j$ of the multi-Harnack smoothing
$C_{\tr_1}$ of $(C, 0)$ is contained in a box parallel to the
coordinate axis and with two vertices of the form
 $(\sim t_j^{e_1 n_1}, \sim t_j^{e_1
m_1})$ and $(\sim t_{j+1}^{e_1 n_1},  \sim t_{j+1}^{e_1 m_1})$,
for $j=1, \dots, g-1$.
\end{Pro}

Now we describe the positions of the ovals in the quadrants of
$\R^2$  of a multi-Harnack smoothing $C_{\tr_1}$ of a real plane
branch $(C,0)$.

\begin{Defi}
We say that the signed topological type of a multi-Harnack
smoothing $C_{\tr_1}$ of a real plane branch $(C,0)$ is normalized
if the signed topological type of the chart of
$P^{\hat{\D}_{g}}_{t_g =1}$ is normalized (see Definition
\ref{normalized-smoothing}).
\end{Defi}

\begin{Pro} \label{normalized-multi}
Let  $C_{\tr_1}$  be  a multi-Harnack smoothing of a real plane
branch $(C,0)$. If the signed topological type of the chart of
$P_{\D_{g}}$ is normalized, then  the same happens for the signed
topological type of the chart $P^{\hat{\D}_j}_{t_j =1}$, for $j=1,
\dots, g-1$.
\end{Pro}
{\em Proof.} The proof is by induction on $g$, using Proposition
\ref{inter} and that $\p_1 (\R^2_{0,0}) \subset \R^2_{0,0}$.
$\Box$

We assume that the signed type of the smoothing is normalized. By
the proof of Theorem \ref{main-top} and the definition an oval of
depth $j$ of $C_{\tr_1}$  arises as a connected component
$\O_{r,s}^{(j)}$ of the chart $\mbox{Ch}_{\D_j}
(P^{\hat{\D}_j}_{t_j =1})$, for $(r,s) \in \mbox{\rm int} (\D_j)
\cap \Z^2$ (see Remark \ref{top-harnack}). An arc of a mixed oval
of depth $j$ corresponds to  a connected component
$\O_{r,s}^{(j)}$ of the chart $\mbox{Ch}_{\D_j}
(P^{\hat{\D}_j}_{t_j =1})$, for $(r,s) \in \Gamma_j \cap \Z^2$.
Since at each stage the smoothing $C^{(j)}_{\tr_j}$ of $C^{(j)}$
is multi-Harnack only the non compact component of the smoothing
$C^{(j)}_{\tr_j}$ meets the coordinate axis. Then we have the
following:

\begin{Pro}
If $C_{\tr_1}$ is a normalized multi-Harnack smoothing of a real
plane branch $(C,0)$ then we can label the ovals of depth $j$
(resp. the mixed ovals of depth $j$)  by $O_{r,s}^{(j)}$ for
$(r,s) \in \mbox{\rm int} (\D_j) \cap \Z^2$  (resp. for $(r,s) \in
\Gamma_j \cap \Z^2$). Then we have that $ O_{r,s}^{(j)} \subset
\R^2_{k,l}$ where
\[
(k, l) = \left\{
\begin{array}{ccc}
 ( e_0 +s, r)  & \mbox{ if } &  j =1
 \\
 (e_{j-1} +s) n_1 n_2\cdots n_{j-1},   (e_{j-1} +s) m_1 n_2\cdots n_{j-1})  & if  & 1 < j \leq
 g.
\end{array}
\right.
\]
\end{Pro}
{\em Proof.} By hypothesis and Proposition \ref{normalized-multi}
we have that the signed topological type of the chart of
$P^{\hat{\D}_j}_{t_j =1}$ is normalized,  for $j=1, \dots, g$. By
Remark \ref{top-harnack2} we have that the component
$\O_{r,s}^{(j)}$ of the chart $\mbox{\rm Ch}^*_{\D_j}
(P^{\hat{\D}_j}_{t_j =1}) $ is contained in the quadrant
$\R^2_{e_{j-1} + s, r}$ with respect to the coordinates $(x_j,
y_j)$. We abuse of notation and denote in the same way the
component $\O_{r,s}^{(j)}$ of the chart $\mbox{\rm Ch}^*_{\D_j}
(P^{\hat{\D}_j}_{t_j =1}) $ and the corresponding connected
component of the smoothing $C_{\tr_j}^{(j)} \cap (\R^*)^2$. By
construction the oval $O_{r,s}^{(j)}$, which appears as a slight
deformation of $\p_1 \circ \cdots \circ \p_{j-1} ( O_{r,s}^{(j)})
$ of $\O_{r,s}^{(j)}$, is contained in the same quadrant $\R^2_{k,
l}$ as  $\p_1 \circ \cdots \circ \p_{j-1} ( \O_{r,s}^{(j)})$.
Notice that on the oval $\p_i \circ \cdots \circ \p_{j-1} (
\O_{r,s}^{(j)}) $ the function $u_i \sim 1$, for $0 < t_i \ll
\cdots \ll t_g \ll 1$ by (\ref{xi-rel}), for $1 \leq i \leq j$.
Then the assertion follows by the definition of the toric maps
$\p_i$ in (\ref{mm-i}). $\Box$

\subsection{Examples}

\begin{figure}[htbp]
 $$\epsfig{file=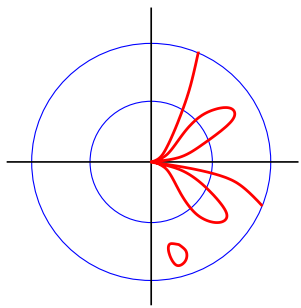, height= 40 mm}$$
 \caption{The deformation $C_{\tr_2}$ inside the Milnor ball \label{ejemplo-def}}
 \end{figure}

\begin{figure}[htbp]
 $$\epsfig{file=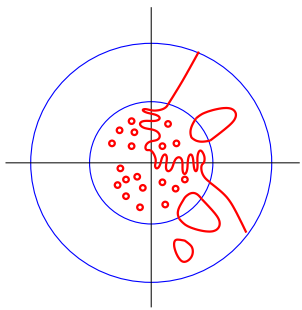, height= 40 mm}$$
 \caption{A multi-Harnack smoothing \label{ejemplo-def-2}}
 \end{figure}

\begin{Exam}
We consider first the constructions of a multi-Harnack smoothing
of the real plane branch $(C, 0)$ defined by the polynomial $F=
(y^2- x^3)^3 - x^{10}$ studied in the Example \ref{10} and
\ref{11}. The Milnor number is equal to $44$.
\end{Exam}
The strict transform $C^{(2)}$ of $(C,0)$ is a simple cusp. Then a
smoothing  $C^{(2)}_{\tr_2}$ of normalized type is of the form
indicated in Figure \ref{hk}.  We indicate the form of the
deformation $C_{\tr_2}$ inside a Milnor ball $B$ for $(C,0)$ in
Figure \ref{ejemplo-def}; the small circle denotes a Milnor ball
$B'$ for the singularity obtained. The smoothing $C_{\tr_1}$ is
the result of perturbing this singularity inside its Milnor ball
and appears in Figure \ref{ejemplo-def-2}. Notice that the ovals
which appear in the smaller ball are infinitesimally smaller than
the others. In this example there is only one oval of depth $2$
and two mixed ovals of depth $1$.

\begin{Exam}
We consider first the constructions of multi-Harnack smoothings of the real plane branch $(C, 0)$
defined by the polynomial $F= (y^2- x^3)^4 - x^{12}y$.  The Milnor number is equal to $86$.
\end{Exam}

\begin{figure}[htbp]
 $$\epsfig{file=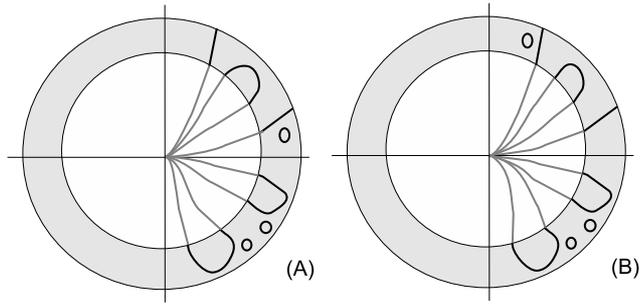, height= 40 mm}$$
 \caption{ Figure (A) corresponds to the normalized case \label{normalized2}}
 \end{figure}

By computing as in Section \ref{toric} we find that   $g=2$,
$(n_1, m_1) = (2, 3)$ and $(n_2, m_2) = (4, 3)$.     It follows
that the strict transform  $ C^{(2)}$ of $(C, 0)$ has local Newton
polygon with vertices $(0,4)$ and $(3,0)$. Figure
\ref{normalizado} shows the signed topological types of a Harnack
smoothing $C^{(2)}_{\tr_2}$  of $(C^{(2)}, o_2)$. Then, the
deformation $C_{\tr_2}$, obtained by Proposition \ref{inter} is
shown in Figure \ref{normalized2} (A) and (B) (notice that all the
branches of the singularity at the origin have the same tangent
line, the horizontal axis, though this is not represented in the
Figure).

\begin{figure}[htbp]
 $$\epsfig{file=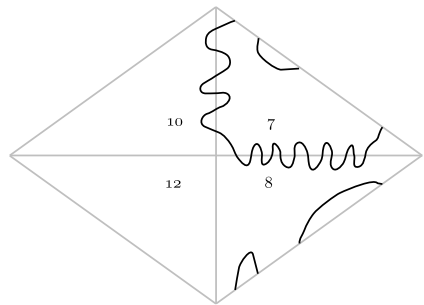, height= 40 mm}$$
 \caption{ The numbers  indicate the   ovals to be added in the corresponding region \label{chart-example}}
 \end{figure}

\begin{figure}[htbp]
 $$\epsfig{file=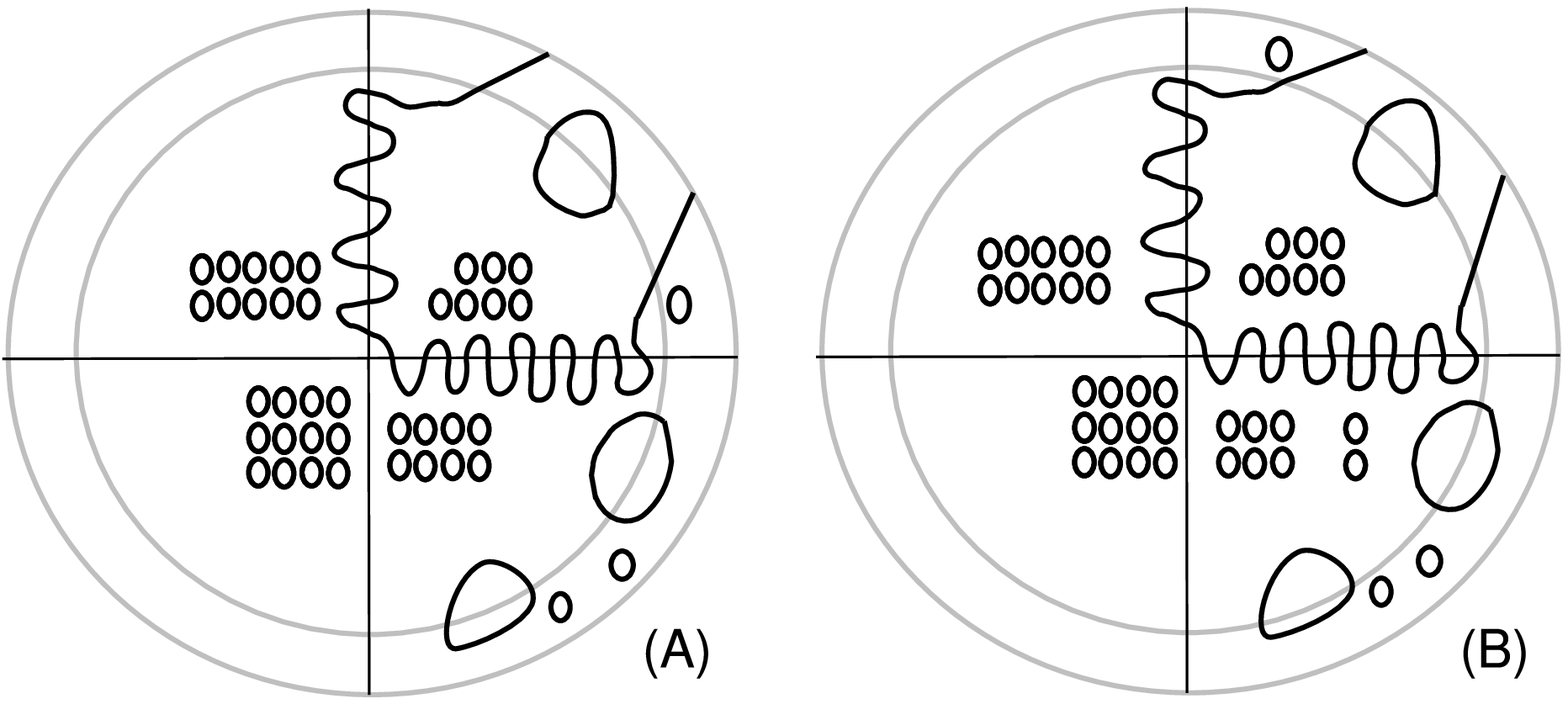, height= 40 mm}$$
 \caption{   \label{normalizado3}}
 \end{figure}

Then, if $C_{\tr_2}$ defines a multi-Harnack smoothing of $(C,0)$
we have that $P^{\hat{\D}_1}_{t_1 =1} (x, y)$ has the chart
represented in Figure \ref{chart-example}. The topology of the
multi-Harnack smoothing is shown in Figure \ref{normalizado3}.
Notice that, as stated in Theorem \ref{main-top}, the topological
type of the resulting msqh-smoothing  $C_{\tr_2}$ is the same in
cases (A) and (B) (meanwhile the signed topological types are
different). The ovals inside the fist ball are of depth $1$, those
intersecting the boundary are mixed ovals of depth $1$ and the
ovals in between both balls are of depth $2$.

\begin{Exam} \label{nested}
We show a Harnack smoothing of the real plane branch $(C, 0)$
defined by the polynomial $F= (y^2- x^3)^7 - x^{24}$ which is not
multi-Harnack.
\end{Exam}

By computing as in Section \ref{toric} we find that $g=2$, $(n_1,
m_1) = (2, 3)$ and $(n_2, m_2) = (7, 6)$ and $\mu(C)_0 = 296$. It
follows that the strict transform  $ C^{(2)}$ of $(C, 0)$ has
local Newton polygon with vertices $(0,7)$ and $(6,0)$.

\begin{figure}[htbp]
 $$\epsfig{file=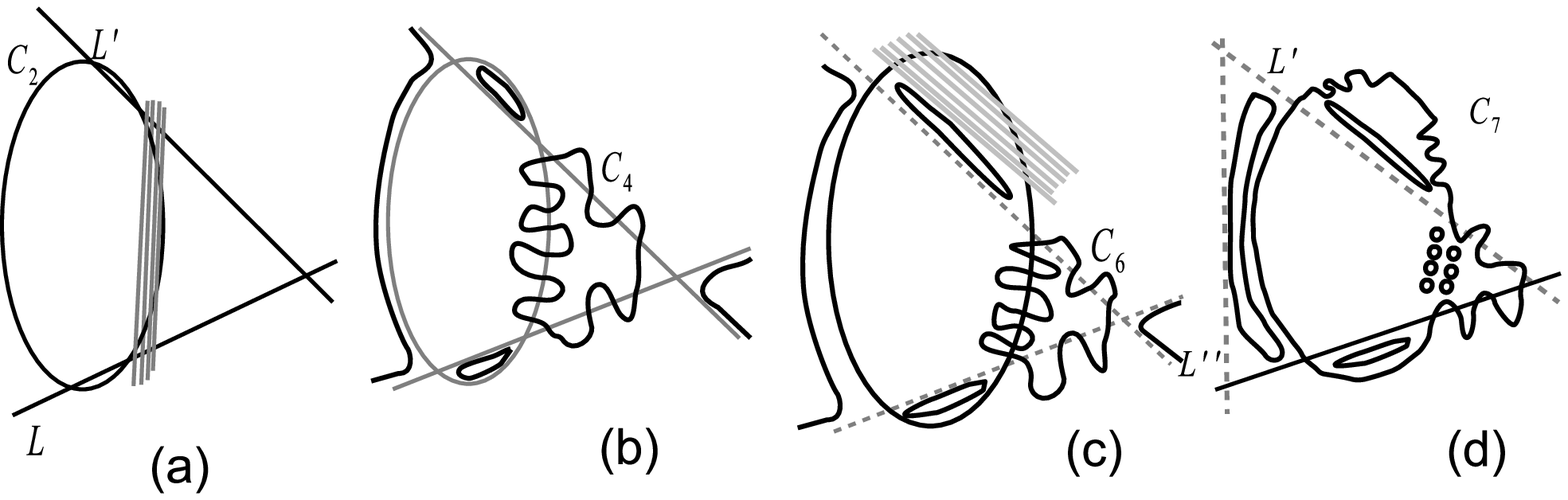, height= 40 mm}$$
\caption{ \label{paso}}
 \end{figure}

\begin{figure}[htbp]
 $$\epsfig{file=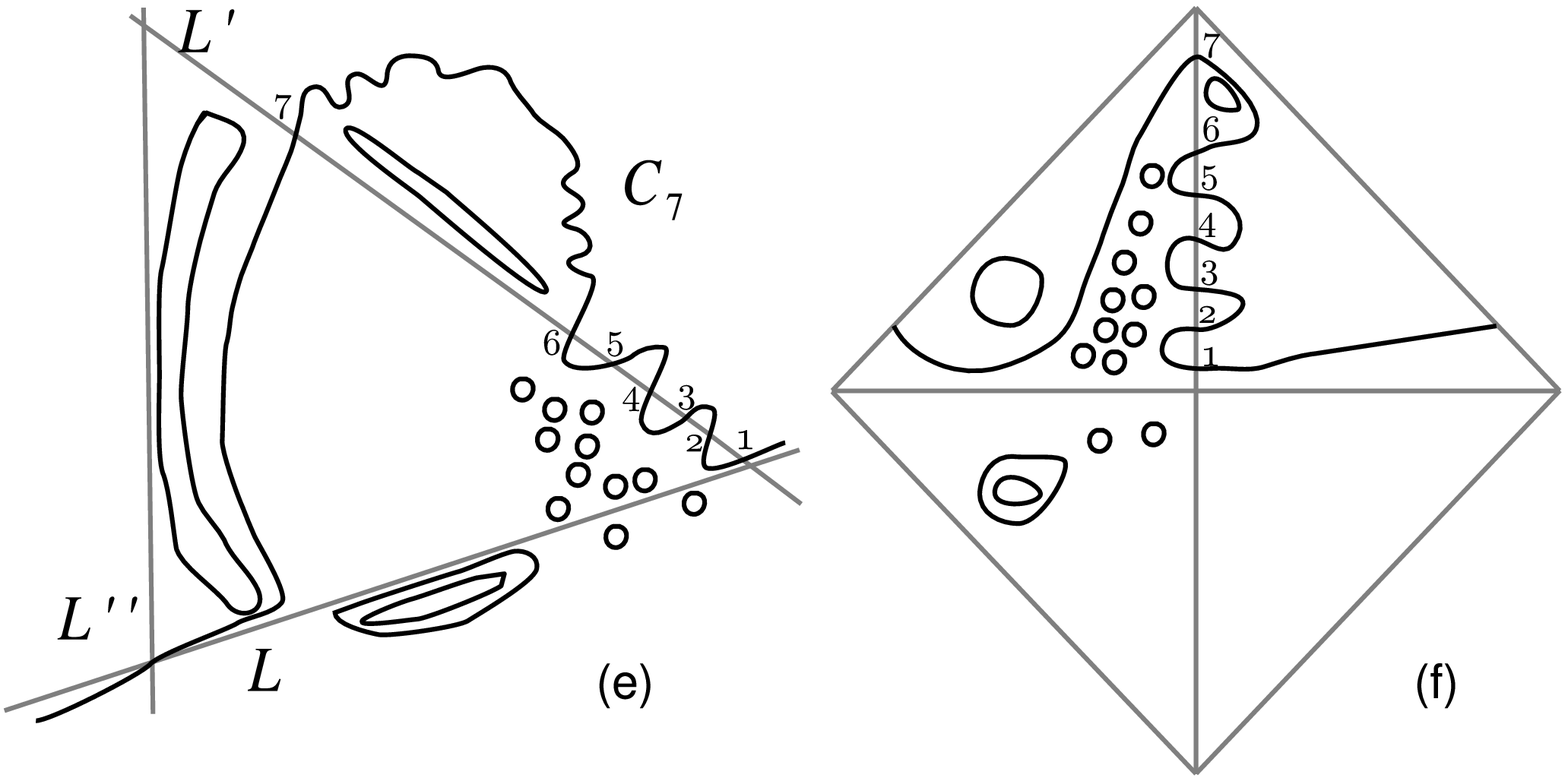, height= 40 mm}$$
\caption{  \label{paso-e1}}
 \end{figure}

We exhibit first a smoothing     $C^{(2)}_{\tr_2}$, defined by a
degree seven curve with Newton polygon equal to $\D_2$. The
construction begins
 by perturbing a degree four
curve, composed of a smooth conic $C_2$  and two lines $L$ and
$L'$  in Figure \ref{paso} (a), with four lines, shown in grey,
intersecting the conic in two real points (see \cite{VirL} for a
summary of construction of real curves by small perturbations).
The result is a smooth quartic $C_4$, as   in Figure  \ref{paso}
(b), where we have indicated in gray the reference lines $L$ and
$L'$ with the conic $C_2$. Then we perturb the union of $C_2$ and
$C_4$, by taking six lines, as in Figure \ref{paso} (c). The
result is a $M$-sextic in maximal position with respect to the
line $L$. The union of both curves is a degree seven curve $C_7$,
as shown
  in Figure  \ref{paso}(d),
where we indicate also the reference lines $L'$ and $L''$. Notice
that
                      in   Figures \ref{paso} (c) and   (d).
the line at infinity changes. See also the construction of curves
corresponding to Figure 13 and 14 of \cite{VirL}. The result of a
suitable perturbation of the singularities of $C_7$ is the
$M$-degree seven curve shown in Figure \ref{paso-e1} (e). Notice
that this curve, which we call also $C_7$, is in maximal position
with respect to the line $L'$ and has maximal intersection
multiplicity with the line $L$ at the point $L \cap L''$. It
follows that a polynomial defining $C_7$, with respect to affine
coordinates $(x, y)$ such that $x=0$ defines the line $L$, $y=0$
defines $L'$ and $L''$ is the line at infinity, has generically
Newton polygon with vertices $(0,0)$, $(7,0)$ and $(0,6)$, since
$C_7$ passes by the point $L \cap L''$. The chart of $C_7$, with
respect to its Newton polygon, is shown in Figure \ref{paso-e1}
(f).

\begin{figure}[htbp]
 $$\epsfig{file=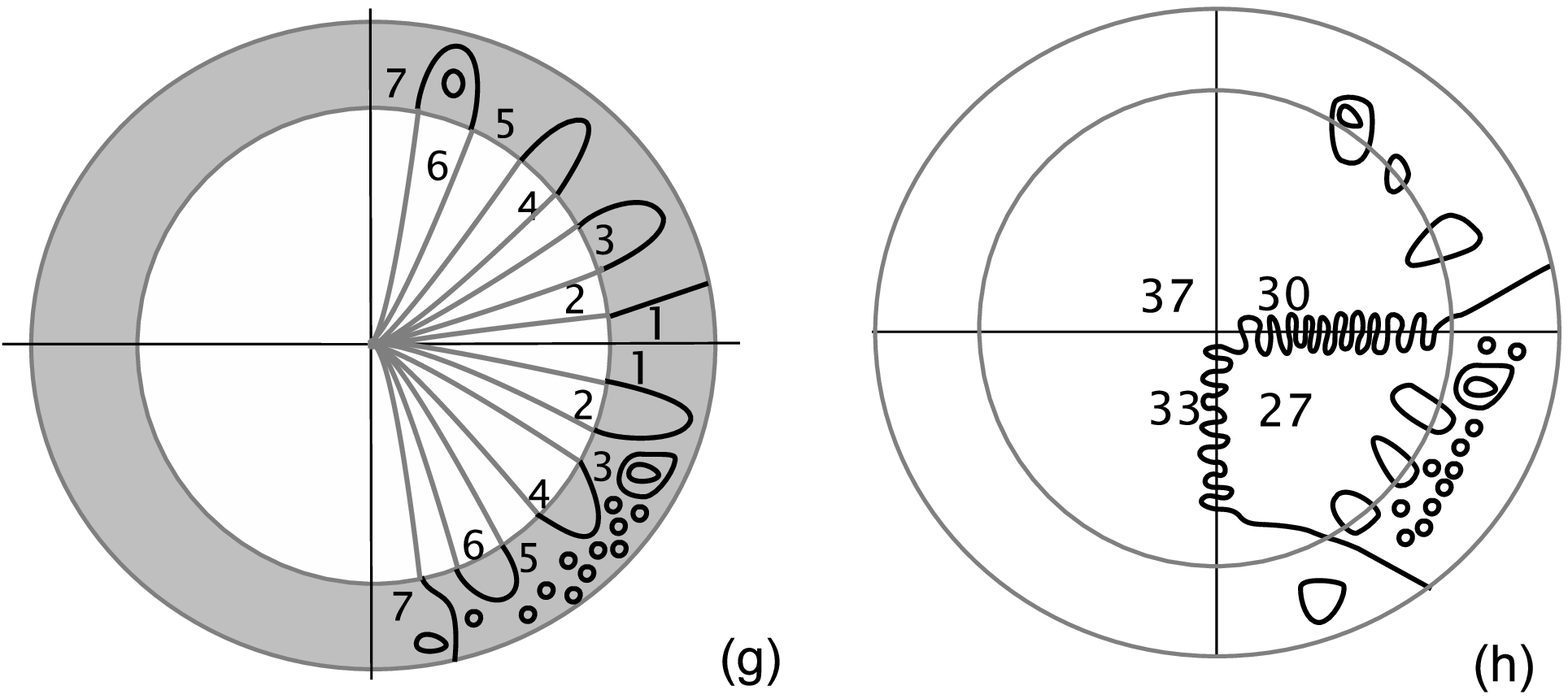, height= 40 mm}$$
\caption{ \label{harnack-non-multi}}
 \end{figure}

We construct from the curve $C_7$  a sqh-smoothing $C^{(2)}_{
\tr_2}$ of $(C^{(2)}, o_2)$. By Proposition \ref{inter} the
topology of the deformation $C_{ \tr_2}$ can be seen from the
chart associated to $C_{\tr_2}$ (compare Figures \ref{paso-e1} (f)
and \ref{harnack-non-multi} (g), where we have indicated by the
same number the corresponding peripheral roots)
 We define then a  msqh-Harnack smoothing $C_{\tr_1}$ of $C_{\tr_2}$
by gluing together the chart of $C_{\tr_2}$ with the chart of a
suitable Harnack smoothing of  $(C_{\tr_2},0)$ (see Theorem
\ref{1par}) . The topology of the msqh-Harnack smoothing
$C_{\tr_1}$ is shown in Figure \ref{harnack-non-multi} (h).


\newpage

 {\bf Acknowledgement.} The authors are grateful to Erwan
 Brugall\'e for suggesting example
 \ref{nested}.

 {\small

\end{document}